\documentclass[11pt]{article}

\usepackage{url}
\usepackage{hyperref}
\usepackage{fullpage}
\usepackage[english]{babel}
\usepackage{amsfonts}
\usepackage{amssymb}
\usepackage{amsmath}
\usepackage{amsthm}
\usepackage{amsopn}
\usepackage{graphicx}
\usepackage{times}

\newcommand{\finishline}{\hfill\hbox{}\linebreak[4]}
\newcommand{\E}{\mathbb E}
\newcommand{\N}{\mathbb N}
\newcommand{\Z}{\mathbb Z}

\newcommand{\R}{\mathbb R}

\newcommand{\GL}{\mathsf{GL}}
\newcommand{\gldz}{\GL_d(\Z)}
\newcommand{\gldr}{\GL_d(\R)}
\newcommand{\odr}{\mathsf{O}_d(\R)}
\newcommand{\sd}{{\cal S}^d}
\newcommand{\sdo}{{\cal S}^d_{>0}}
\newcommand{\sdgeo}{{\cal S}^d_{\geq 0}}
\newcommand{\SC}{{\mathbf\Delta}}
 
\newcommand{\MA}{{\cal A}}
\newcommand{\MC}{{\cal C}}
\newcommand{\MD}{{\cal D}}

\newcommand{\MN}{{\cal N}}

\newcommand{\MS}{{\cal S}}
\newcommand{\MR}{{\cal R}}
\newcommand{\MT}{{\cal T}}

\DeclareMathOperator{\trace}{trace}
\DeclareMathOperator{\vol}{vol}
\DeclareMathOperator{\conv}{conv}
\DeclareMathOperator{\cone}{cone}

\DeclareMathOperator{\lin}{lin}
\DeclareMathOperator{\hess}{hess}
\DeclareMathOperator{\gradient}{grad}
\DeclareMathOperator{\bd}{bd}

\DeclareMathOperator{\dist}{dist}
\DeclareMathOperator{\BR}{BR}
\DeclareMathOperator{\CM}{CM}
\DeclareMathOperator{\Del}{Del}
\DeclareMathOperator{\vertex}{vert}

\DeclareMathOperator{\DV}{DV}
\DeclareMathOperator{\Aut}{Aut}
\DeclareMathOperator{\opt}{opt}

\theoremstyle{definition}
\newtheorem{definition}{Definition}[section]
\newtheorem{proposition}[definition]{Proposition}

\newtheorem{theorem}[definition]{Theorem}
\newtheorem{corollary}[definition]{Corollary}
\newtheorem{conjecture}[definition]{Conjecture}
\newtheorem{example}[definition]{Example}
\newtheorem{question}[definition]{Question}
\newtheorem{problem}[definition]{Problem}

\renewcommand{\vec}[1]{\boldsymbol{#1}}

\author{Achill Sch\"urmann and Frank Vallentin\footnote{This work was partially
supported by the Edmund Landau Center for Research in Mathematical
Analysis and Related Areas, sponsored by the Minerva Foundation
(Germany), and by the ``Deutsche Forschungsgemeinschaft'' (DFG) under grant SCHU 1503/4-1.}}

\title{
Computational Approaches\\
to Lattice Packing and Covering Problems}

\date{27th July 2005}

\begin{document}

\selectlanguage{english}

\maketitle

\begin{abstract}
We describe algorithms which address two classical problems in lattice
geometry: the lattice covering and the simultaneous lattice
packing-covering problem. Theoretically our algorithms solve the two
problems in any fixed dimension $d$ in the sense that they approximate
optimal covering lattices and optimal packing-covering lattices within
any desired accuracy. Both algorithms involve semidefinite
programming and are based on Voronoi's reduction theory for positive
definite quadratic forms, which describes all possible Delone
triangulations of~$\Z^d$.

In practice, our implementations reproduce known results in
dimensions~$d \leq 5$ and in particular solve the two problems in
these dimensions. For~$d = 6$ our computations produce new best known
covering as well as packing-covering lattices, which are closely
related to the lattice~$\mathsf{E}^*_6$. For~$d = 7,8$ our approach
leads to new best known covering lattices. Although we use numerical
methods, we made some effort to transform numerical evidences into
rigorous proofs. We provide rigorous error bounds and prove that some
of the new lattices are locally optimal.
\end{abstract}




\newpage

\tableofcontents

\newpage


\section{Overview}
\label{sec:overview}

Two classical problems in the geometry of numbers are the
determination of the most economical lattice sphere packings and
coverings of the Euclidean $d$-space~$\E^d$. In this paper we describe
algorithms for the lattice covering and the simultaneous lattice
packing and covering problem (lattice packing-covering problem in the
sequel).

Roughly speaking, both problems are concerned with the most economical
way to cover~$\E^d$.  In the case of the lattice covering 
problem, the
goal is to maximize the volume of a fundamental domain in a lattice
covering with unit spheres. Roughly speaking, we want to minimize the
number of unit spheres which are needed to cover arbitrarily large but
finite regions of $\E^d$.  The objective of the lattice
packing-covering problem is to maximize the minimal distance between
lattice points in a lattice covering with unit spheres.  

The aim of this paper is to give an introduction to the mathematical
tools that allow us, at least in theory, to solve the two problems
computationally. For a fixed dimension $d$, our algorithms approximate
optimal covering lattices and optimal packing-covering lattices within
any desired accuracy. In this overview we want to describe the
structure of the paper.

The basic concepts and notations, which we use throughout the paper,
are in Section~\ref{sec:notations}. There we also give a precise
definition of the two problems under consideration. The reader familiar with
sphere packings and coverings, as well as lattices and their relation
to positive definite quadratic forms, may skip this section.

In Section~\ref{sec:lcp_overview} we review known results and the
history of the lattice covering problem, and in
Section~\ref{sec:lspc_overview} we review known results and the
history of the packing-covering problem.

Our algorithms as well as the known results by other authors are
mainly based on a reduction theory for positive definite quadratic
forms by Voronoi. We give a detailed description of this main
ingredient in Section~\ref{sec:voronoireduction} with a special focus
on computational implementability.

The other main tool comes from convex optimization
theory. Semidefinite programming problems and determinant maximization
problems are briefly described in Section~\ref{sec:convopt}. We
describe how duality theory together with rational approximations can
be used to provide rigorous error bounds. Both problems have in common
that one has to minimize a convex function on variables that satisfy
some linear matrix inequalities (LMIs).

In Section~\ref{sec:lmi} we describe how the constraint that a lattice
gives a \textit{unit} sphere covering can be modeled by LMIs.  

In Section~\ref{sec:algorithms} we combine these tools and attain
algorithms which theoretically solve the two problems under
consideration.

Due to a combinatorial explosion of the number of different Delone
triangulations, our implementations of the algorithms only give
complete solutions for~$d \leq 5$.  Moreover, the convex optimization
algorithms we used are interior point methods and so yield only
approximations. Therefore, in Section~\ref{sec:classical}, we collect
some mathematical tools which allow us to determine exact results from
these approximations.  In particular we can test computationally
whether or not a given positive definite quadratic form gives a
locally optimal solution. In the case of the lattice packing-covering
problem, we can test if such a solution is isolated.

In order to run a heuristic search for good lattices, it is necessary
to have local lower bounds that we can compute fast. We describe one
class of such bounds depending on the methods of inertia in
Section~\ref{sec:moments}, which we used to find new lattices in
dimension~$6$.

Both problems have been previously solved only for
dimensions~$d \leq 5$.  Our implementations not
only verify all of these results, but also
attain additional information on locally optimal solutions for $d
= 5$. Moreover, we find new best known
lattices for both problems in dimension~$6$, $7$, and $8$.  In
particular, we answer an open question by Ryshkov.  In
Section~\ref{sec:compute} we report on our results.  There we
distinguish between conjectures for which we only have numerical
evidence and theorems for which we have rigorous proofs.

\newpage


\section{Basic Concepts and Notations}
\label{sec:notations}

In this section we fix the notation we use throughout this paper.
We refer the reader to \cite{br-1979}, \cite{gl-1987} and
\cite{cs-1988} for further information about the introduced concepts.
The reader familiar with sphere packings and coverings, as well as
lattices and their relation to positive definite quadratic forms, may
skip this section.

Let $\E^d$ be a $d$-dimensional Euclidean space equipped with inner
product $\langle \cdot, \cdot \rangle$, norm $\|\cdot\|$ and unit ball
$B^d=\left\{ \vec{x} \in \E^d : \|\vec{x}\| \leq 1 \right\}$.  A
\textit{lattice} $L$ is a discrete subgroup in $\E^d$.  From now on
we assume that all lattices $L$ have full rank $d$; that is,
there exists a regular matrix $A \in \gldr$ with $L = A\Z^d$.  The
columns of the matrix~$A$ are called a \textit{basis} of~$L$.  All
bases of $L$ are of the form $AU$ with $U \in \gldz$.  Thus, the
\textit{determinant} $\det(L) = |\det(A)| > 0$ of the lattice $L$ is
well defined. We say that two lattices $L$~and~$L'$ belong to the same
\textit{isometry class} if, for every basis~$A$ 
of~$L$, there is a
basis~$A'$ of~$L'$ and an orthogonal transformation $O \in \odr$ 
such that $A' = OA$.

The Minkowski sum $L + \alpha B^d = \{\vec{v} + \alpha\vec{x} :
\vec{v} \in L, \vec{x} \in B^d\}$, with $\alpha \in \R_{>0}$, is a
\textit{lattice packing} if the translates of $\alpha B^d$ have
mutually disjoint interiors and a \textit{lattice covering} if $\E^d =
L + \alpha B^d$. The \textit{packing radius} $\lambda(L)$ of a lattice
$L$ is given by
\[
\lambda(L)
=
\max\{ \lambda : \mbox{$L + \lambda B^d$ is a lattice packing} \},
\]
and the \textit{covering radius} $\mu(L)$ by
\[
\mu(L)
=
\min \{ \mu : \mbox{$L + \mu B^d$ is a lattice covering} \}.
\]
The above values are attained: The packing radius is equal to half the
length of a shortest non-zero vector of~$L$ and the covering radius is
equal to the maximum distance of points in $\E^d$ to a closest lattice
vector. The packing radius is the inradius of the
\textit{Dirichlet-Voronoi polytope} of $L$
\[
\DV(L) = 
\{ \vec{x} \in\E^d : \|\vec{x}\| \leq \|\vec{x}-\vec{v}\|, \vec{v} \in L \},
\]
and the covering radius is its circumradius. Both functionals are
homogeneous, that is, for $\alpha \in \R$ we have
\[
\mu(\alpha L) = |\alpha|\mu(L)
\quad\mbox{and}\quad
\lambda(\alpha L) = |\alpha|\lambda(L)
.
\]
Thus, the \textit{covering density} 
\[
\Theta(L)=\frac{\mu(L)^d}{\det(L)}\cdot\kappa_d
,\quad \kappa_d = \vol B^d,
\]
is invariant with respect to scaling of~$L$.  The same is true for the
\textit{packing-covering constant}
\[
\gamma(L)=\frac{\mu(L)}{\lambda(L)}.
\]
All these functionals are invariants of the isometry classes. 
 
\noindent 
In this paper we study the following two problems:

\begin{center}
\fbox{
\begin{minipage}{15.5cm}
\begin{problem} (\textit{Lattice covering problem})
\label{problem:cov}
\finishline
For a given $d\in\N$, determine
$\Theta_d = \min_L \Theta(L)$,
where $L \subseteq \E^d$ runs over all $d$-dimensional lattices. 
\end{problem}

\begin{problem} (\textit{Lattice packing-covering problem})
\label{problem:packcov}
\finishline
For a given $d\in\N$, determine $\gamma_d = \min_L \gamma(L)$, where
$L \subseteq \E^d$ runs over all $d$-dimensional lattices.
\end{problem}
\end{minipage}
}
\end{center}

We describe history and results of both problems in
Section~\ref{sec:lcp_overview} and Section~\ref{sec:lspc_overview}.
There we assume that the reader is familiar with certain important
lattices which are described in \cite{cs-1988}, Ch.~4. Historically
and for computational reasons, the problems were studied by using the
intimate relation between lattices and positive definite quadratic
forms (PQFs).
 
Let us describe this relation: To a $d$-dimensional lattice $L=A\Z^d$
with basis $A$ we associate a $d$-dimensional PQF
\[
Q[\vec{x}]
=\vec{x}^t A^t A \vec{x}
=\vec{x}^t G \vec{x}
,
\]
where the Gram matrix $G = A^t A$ is symmetric and positive definite.
We abuse notation and identify quadratic forms with symmetric matrices by
saying $Q = G$ and $Q[\vec{x}] = \vec{x}^t Q \vec{x}$.  The set of
quadratic forms is a $\frac{d(d+1)}{2}$ dimensional real vector space
$\sd$, in which the set of PQFs forms an open, convex cone $\sdo$.
Its closure is the convex cone of all positive semidefinite
quadratic forms $\sdgeo$, which is pointed at $\vec{0}$.

Note that $Q$ depends on the chosen basis $A$ of $L$.  Two arbitrary
bases $A$ and $B$ of $L$ are transformed into each other by a
unimodular transformation, that is, there exists an $U\in\gldz$ such
that $A=BU$.  Thus, $\gldz$ acts on $\sdo$ by $ Q \mapsto U^t Q
U$. Two PQFs lying in the same orbit under this action are called
\textit{arithmetically equivalent}. This definition naturally extends
to positive semidefinite quadratic forms.

Thus every lattice uniquely determines an arithmetical equivalence
class of PQFs. On the other hand, every PQF $Q$ admits a Cholesky
decomposition $Q=A^t A$, where the upper triangular matrix $A$ is
uniquely determined up to an orthogonal transformation $O \in
\odr$. Altogether, we have a bijection between isometry classes of
lattices $\odr \backslash \gldr / \gldz$ and arithmetical equivalence
classes of PQFs $\sdo /\gldz$.

As a consequence, the lattice covering and the lattice
packing-covering problem translate into problems for PQFs: The
\textit{determinant} (or discriminant) of a PQF $Q$ is defined by
$\det(Q)$.  The \textit{homogeneous minimum} $\lambda(Q)$ and the
\textit{inhomogeneous minimum} $\mu(Q)$ are given by
\[
\lambda(Q) = \min_{\vec{v} \in \Z^d \backslash \{\vec{0}\}} Q[\vec{v}],
\quad
\mbox{and}
\quad
\mu(Q) = \max_{\vec{x} \in \E^d} \min_{\vec{v} \in \Z^d} Q[\vec{x} - \vec{v}].
\]
A corresponding lattice $L$ satisfies $\det(L) = \sqrt{\det(Q)}$,
$\mu(L)=\sqrt{\mu(Q)}$, $\lambda(L)=\sqrt{\lambda(Q)}/2$.  Therefore
our goal is to minimize
\[
\Theta(Q) = \Theta(L) = \sqrt{\frac{\mu(Q)^d}{\det Q}} \cdot \kappa_d,
\quad
\mbox{and}
\quad
\gamma(Q) = \gamma(L) = 2 \cdot \sqrt{\frac{\mu(Q)}{\lambda(Q)}}
\]
among all PQFs $Q\in\sdo$. 

Since $\Theta(Q)$ and $\gamma(Q)$ are invariant with respect to
the action of $\gldz$ on $\sdo$, we only need to consider one PQF in
each arithmetical equivalence class. Finding a fundamental domain in
$\sdo$ is one of the most basic and classical problems in 
the geometry of numbers.  Such a \textit{reduction theory} for PQFs, especially
suitable for Problem~\ref{problem:cov} and
Problem~\ref{problem:packcov}, is due to Voronoi. We describe it in
detail in Section~\ref{sec:voronoireduction}.


\section{The Lattice Covering Problem}
\label{sec:lcp_overview}

Kershner, in~1939, was the first to consider the lattice convering problem.
In \cite{kershner-1939} he showed that the hexagonal lattice
(see Figure~1) gives the most economical sphere covering in the plane
even without the restriction of being a lattice covering.

\begin{center}
\begin{center}
\fbox{
\includegraphics[width=9cm]{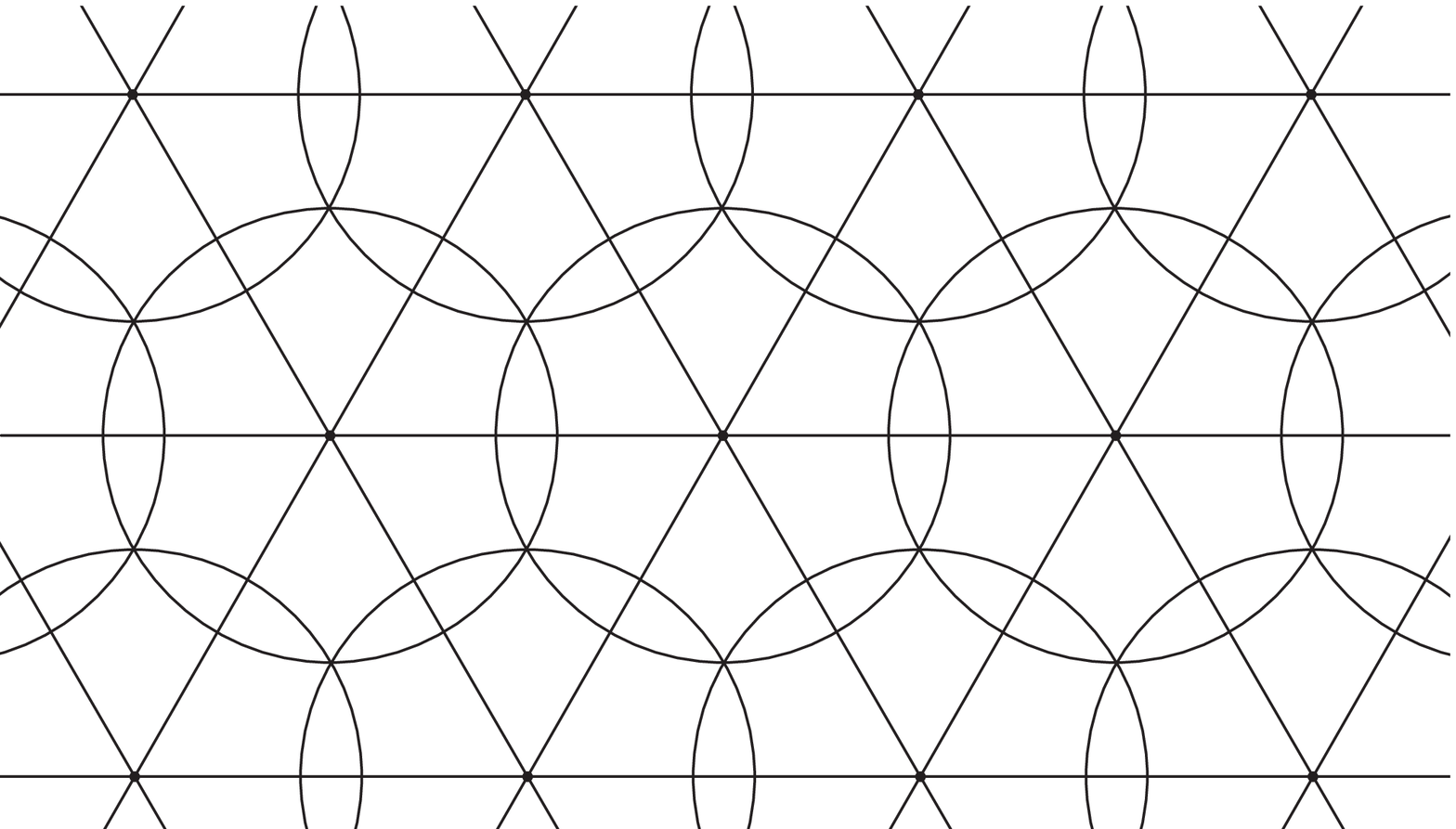}
}
\end{center}
\smallskip
\textbf{\textsf{Figure 1.}} The sphere covering given by the hexagonal lattice.
\end{center}

Since then the lattice covering problem has been solved up to
dimension~$5$ (see Table~1). In all these cases the lattice ${\mathsf
A}_d^*$, whose covering density equals
$$
\Theta({\mathsf A}_d^*) = 
\sqrt{\left(\frac{d(d+2)}{12(d+1)}\right)^d (d+1)} \cdot \kappa_d,
$$ 
provides the optimal lattice covering. Gameckii
\cite{gameckii-1962}, \cite{gameckii-1963}, and Bleicher
\cite{bleicher-1962} were the first to compute the covering density
of ${\mathsf A}_d^*$ for general $d$. They also showed that it is
locally optimal with respect to covering density in every dimension.

The optimality of the body centered cubic lattice~${\mathsf A}_3^*$
whose Dirichlet-Voronoi polytope is a regular truncated octahedron
(the Dirichlet-Voronoi polytope of ${\mathsf A}_d^*$ is a regular
permutahedron) was first proven by Bambah \cite{bambah-1954a}.  Later,
Barnes substantially simplified Bambah's proof in \cite{barnes-1956}
and strengthened the result by showing that in dimensions~$2$ and~$3$
the lattice~${\mathsf A}_d^*$ is the unique locally optimal lattice
covering.  He used Voronoi's reduction theory and anticipated that
this is the right setup for solving the lattice covering problem in
dimensions larger than three.  Our algorithm in
Section~\ref{ssec:solvinglcp} confirms his anticipation.  A third
proof, which is mainly elementary and unlike the previous two does not
use any reduction theory of PQFs, was given by Few \cite{few-1956}. At
the moment no attempt is known to the authors to show that the optimal
three-dimensional lattice covering also gives an optimal sphere
covering without lattice restriction.

In \cite{bambah-1954b} Bambah conjectured that the lattice ${\mathsf
A}_4^*$ gives the least dense four-dimensional lattice covering. In
\cite{dr-1963} Delone and Ryshkov proved Bambah's conjecture. In
\cite{baranovskii-1965}, \cite{baranovskii-1966} Baranovskii gave an
alternative proof of this fact. He determined all locally optimal
lattice coverings in dimension~$4$. Dickson \cite{dickson-1967} gave
another alternative proof of this fact.

In a series of papers \cite{ryshkov-1973a}, \cite{br-1973},
\cite{br-1975}, \cite{rb-1975} Ryshkov and Baranovskii solved the
lattice covering problem in dimension~$5$. They prepared a $140$-page
long monograph \cite{rb-1976} based on their investigations.

\begin{center}
\begin{tabular}{|c|c|c|}
\hline
${\mathbf d}$ & \textbf{lattice covering} & \textbf{covering density $\mathbf \Theta_d$}\\ 
\hline
 $1$ & $\Z^1$               &  $1$        \\
 $2$ & ${\mathsf A}^*_2$    &  $1.209199$ \\
 $3$ & ${\mathsf A}^*_3$    &  $1.463505$ \\
 $4$ & ${\mathsf A}^*_4$    &  $1.765529$ \\
 $5$ & ${\mathsf A}^*_5$    &  $2.124286$ \\
\hline
\end{tabular}
\par\smallskip
\textbf{\textsf{Table 1.}} Optimal lattice coverings.
\end{center}

In \cite{ryshkov-1967} Ryshkov raised the question of finding the
lowest dimension~$d$ for which there is a better lattice covering than
the one given by ${\mathsf A}^*_d$. In the same paper he showed that
${\mathsf A}^*_d$ is not the most efficient lattice covering for all
even $d \geq 114$ and for all odd $d \geq 201$. One of our main 
results in this paper is the answer to Ryshkov's question (see
Theorem~\ref{theorem:ryshkov-answer}): Dimension~$d = 6$ is the lowest
dimension for which there is a better lattice covering than the one
given by ${\mathsf A}^*_d$.  The proof is based on a computer
search. In Section~\ref{ssec:solvinglcp} we will give an algorithm
which finds all locally optimal lattice coverings in a given
dimension. Using this algorithm we are able to verify all of the known
results about optimal lattice coverings up to
dimension~$5$. Unfortunately, due to a combinatorial explosion, the
algorithm can not be applied practically in dimension~$6$ or
greater. Nevertheless, we are able to find good lattice coverings in
dimension~$6$, $7$, and $8$ by applying several heuristics. We will
give more details in Section~\ref{sec:compute}.

\begin{center}
\begin{tabular}{|c|c|c|}
\hline
$\mathbf d$ & \textbf{lattice} & \textbf{covering density $\mathbf\Theta$} \\ 
\hline 
 $6$ & ${\mathsf L}^{c1}_6$ & $2.464801$ \\ 
 $7$ & ${\mathsf L}^{c}_7$  & $2.900024$ \\ 
 $8$ & ${\mathsf L}^{c}_8$  & $3.142202$ \\ 
 $9$ & ${\mathsf A}^5_9$    & $4.340185$ \\ 
$10$ & ${\mathsf A}^*_{10}$ & $5.251713$ \\ 
$11$ & ${\mathsf A}^4_{11}$ & $5.598338$ \\ 
$12$ & ${\mathsf A}^*_{12}$ & $7.510113$ \\ 
$13$ & ${\mathsf A}^7_{13}$ & $7.864060$ \\ 
$14$ & ${\mathsf A}^5_{14}$ & $9.006610$ \\ 
$15$ & ${\mathsf A}^8_{15}$ & $11.601626$ \\ 
$16$ & ${\mathsf A}^*_{16}$ & $15.310927$ \\ 
$17$ & ${\mathsf A}^*_{17}$ & $18.287811$ \\ 
$18$ & ${\mathsf A}^*_{18}$ & $21.840949$ \\ 
$19$ & ${\mathsf A}^*_{19}$ & $26.081820$ \\ 
$20$ & ${\mathsf A}^*_{20}$ & $31.143448$ \\ 
$21$ & ${\mathsf A}^*_{21}$ & $37.184568$ \\ 
$22$ & $\Lambda_{22}^*$     & $\leq 27.8839$ \\ 
$23$ & $\Lambda_{23}^*$     & $\leq 15.3218$ \\ 
$24$ & $\Lambda_{24}$       & $7.903536$ \\
\hline
\end{tabular}
\par\smallskip
\textbf{\textsf{Table 2.}} Least dense known (lattice) coverings
up to dimension $24$.
\end{center}

What else is known? In Table~2 we list all the least dense known
lattice coverings in dimensions~$6$ to~$24$. At the same time this
list gives the least dense known sphere coverings: There is no
covering of equal spheres known which is less dense than the best
known lattice covering. This table is an update of Table~$2.1$ in
\cite{cs-1988} and we provide an up-to-date table on our
web page \cite{sv-2005}.  
We conclude this section by briefly describing the 
origins of these updates.

The Leech lattice $\Lambda_{24}$ yields the best known lattice
covering in dimension $24$.  The covering density of the Leech lattice
was computed by Conway, Parker and Sloane (\cite{cs-1988} Ch.~23). It
is not too brave to conjecture that the Leech lattice gives the
optimal $24$-dimensional sphere covering. In \cite{sv-2004b} we took a
first step towards proving this conjecture by showing that the Leech
lattice gives a locally optimal lattice covering. Using the Leech
lattice Bambah and Sloane constructed in \cite{bs-1982} a series of
lattices in dimensions $d \geq 24$ which give a thinner lattice
covering than ${\mathsf A}_d^*$.  It seems that as a ``corollary'' of
the existence of the Leech lattice the duals of the laminated lattices
$\Lambda_{22}$ and $\Lambda_{23}$ give good lattice coverings. Their
covering densities were estimated by Smith \cite{smith-1988}, but we
do not know the exact values of $\Theta(\Lambda^*_{22})$ and
$\Theta(\Lambda^*_{23})$. For the definitions of the root lattices
$\mathsf{A}_n, \mathsf{D}_n, \mathsf{E}_n$, and the laminated lattices
$\Lambda_n$, we refer to \cite{cs-1988}, Ch.~4 and Ch.~6.

In \cite{coxeter-1951} Coxeter gave a list of locally optimal lattice
\textit{packings} (extreme lattices) which are related to Lie groups.
One of them is the infinite series of locally densest packing lattices
${\mathsf A}^r_d$ where $d \geq 2$ and $r > 1$ divides $d + 1$. The
lattice $\mathsf{A}^r_d$ is the unique sublattice of $\mathsf{A}^*_d$
containing $\mathsf{A}_d$ to index~$r$.  In \cite{baranovskii-1994}
Baranovskii determined the covering density of the lattice covering
given by ${\mathsf A}^5_9$, which is slightly better than the one given
by ${\mathsf A}^*_9$.  Recently, Anzin extended Baranovskii's work. In
\cite{anzin-2002} he computed the covering densities of ${\mathsf
A}^4_{11}$ and ${\mathsf A}^7_{13}$, and in a private communication he
reported on computing the covering densities of ${\mathsf A}^5_{14}$
and ${\mathsf A}^8_{15}$. They all give less dense lattice coverings
than those provided by the corresponding ${\mathsf A}^*_d$.  We do not
know whether these lattice coverings are locally optimal.  

To answer Ryshkov's question exhaustively it will be necessary to
further investigate lattice coverings in the dimensions $d = 10, 12,
16, 17, \ldots, 21$. We hope that the methods we present in this paper
will be useful for this project.


\section{The Lattice Packing-Covering Problem}
\label{sec:lspc_overview}

The lattice packing-covering problem has been studied in different
contexts and there are several different names and interpretations of
the lattice packing-covering constant $\gamma_d$.  Lagarias and
Pleasants \cite{lp-2002} referred to it as the ``Delone
packing-covering constant''. Ryshkov \cite{ryshkov-1974} studied the
equivalent problem of minimizing the density of $(r,R)$-systems. An
$(r, R)$-system is a discrete point set $X \subseteq \R^d$ where (1)
the distance between any two points of $X$ is at least $r$ and (2) the
distance from any point in $\R^d$ to a point in $X$ is at most $R$. If
$X$ is a lattice, then $r = \lambda(X)/2$ and $R = \mu(X)$.

\begin{center}
\begin{center}
\fbox{
\includegraphics[width=9cm]{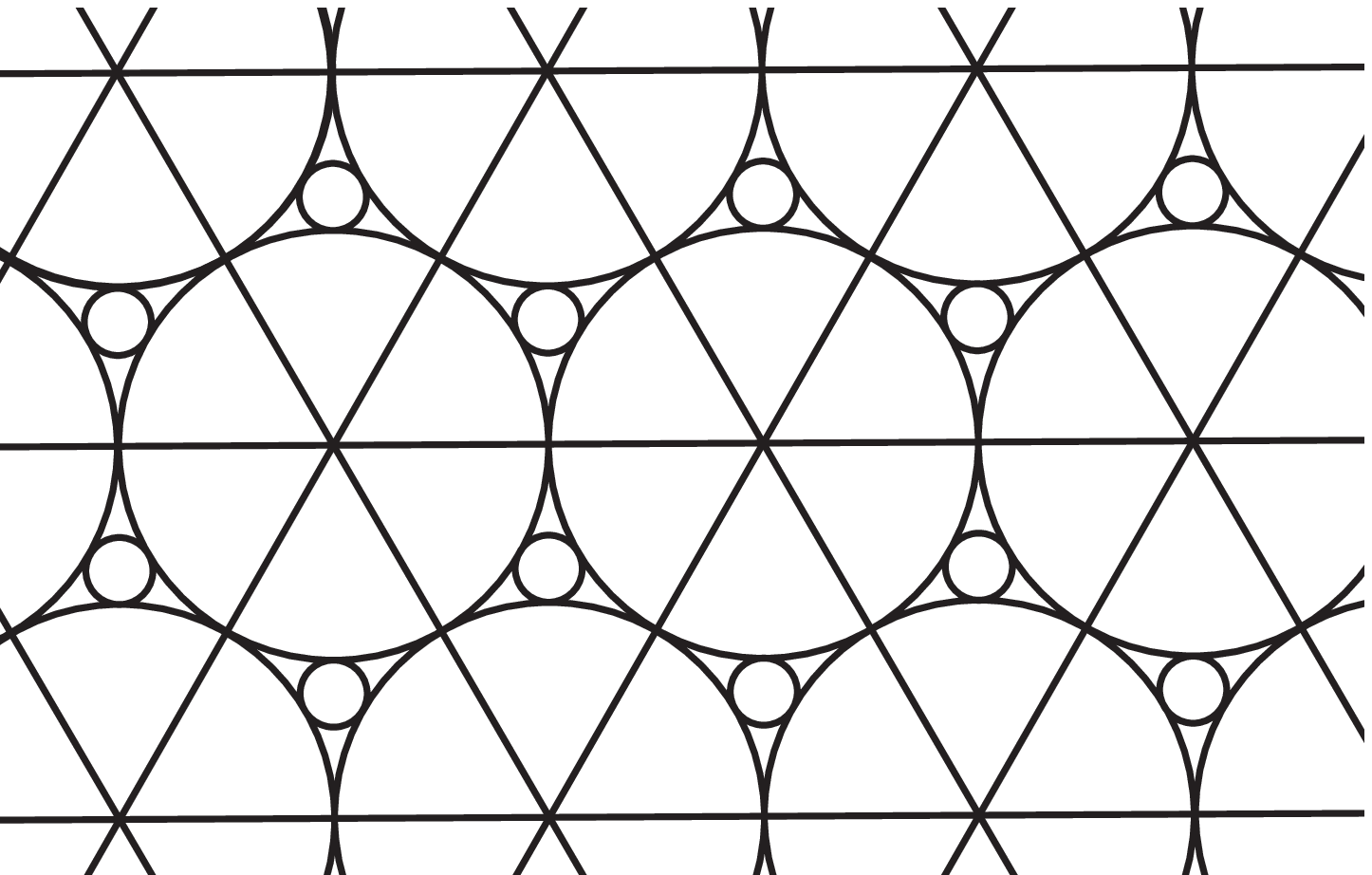}
}
\end{center}
\smallskip
\textbf{\textsf{Figure 2.}} A close sphere packing given by the hexagonal
lattice.
\end{center}

Geometrically we may think of solving the lattice packing-covering problem as
of maximizing the minimum distance between lattice points in a lattice
covering with unit spheres. 
Alternatively, we may think of it as minimizing the
radius of a largest sphere that could additionally be packed into a
lattice packing of unit spheres (see Figure~2).  
This minimal
gap-radius is equal to $\gamma_d - 1$. Therefore the problem raised by
L.\ Fejes T\'oth \cite{fejes-toth-1976}  
of finding ``close packings''
attaining this gap-radius is another formulation of the
packing-covering problem.

The last interpretation shows that $\gamma_d \geq 2$ would imply that
in any $d$-dimensional lattice packing with spheres of unit radius
there is still space for spheres of radius $1$.  In particular, this
would prove that densest sphere packings in dimension~$d$ are
non-lattice packings. This phenomenon is likely to be true for
large dimensions, but it has not been verified for any $d$ so far.

\begin{problem}
\label{problem:pack-cov}
Does there exist a $d$ such that $\gamma_d\geq 2$?
\end{problem}

Note that this problem is particularly challenging in view of the 
asymptotic bound $\gamma_d \leq 2 + o(1)$ due to Butler \cite{butler-1972}.

As for the lattice covering problem, the lattice packing-covering
problem has been solved up to dimension~$5$ (see Table~3).  Ryshkov
\cite{ryshkov-1974} solved the general $2$-dimensional case. The
$3$-dimensional case was settled by B\"or\"oczky
\cite{boeroeczky-1986}, even without the restriction to lattices.  The
$4$- and $5$-dimensional cases were solved by Horv\'ath
\cite{horvath-1980}, \cite{horvath-1986}.  Note that the lattices
$\mathsf{Ho}_4$ and $\mathsf{Ho}_5$ (see Section~\ref{sec:compute})
discovered by Horv\'ath are neither best covering, nor best packing
lattices.  As in the case of lattice coverings, the results were
attained by using Voronoi's reduction theory.

\begin{center}
\begin{tabular}{|c|c|c|}
\hline
$\mathbf d$ & \textbf{lattice} & \textbf{lattice packing-covering constant} 
$\mathbf \gamma_d$ \\ 
\hline
 $1$ & $\Z^1$ & $1$ \\
 $2$ & $\mathsf{A}_2^\ast$  & $\frac{2}{\sqrt{3}} \approx 1.154700$ \\
 $3$ & $\mathsf{A}_3^\ast$  & $\sqrt{5/3} \approx 1.290994$ \\
 $4$ & $\mathsf{Ho}_4$      & $\sqrt{2\sqrt{3}}(\sqrt{3}-1) \approx 1.362500$ \\
 $5$ & $\mathsf{Ho}_5$      & $\sqrt{3/2+\sqrt{13}/6} \approx 1.449456$ \\
\hline
\end{tabular}
\par\smallskip
\textbf{\textsf{Table 3.}} Optimal packing-covering lattices.
\end{center}

Our computations, described in 
Section~\ref{sec:compute}, verify all of
the known results in dimension $\leq 5$.  As for the covering problem,
none of the values $\gamma_d$ has been determined in a dimension
$d\geq 6$ so far.  In Section~\ref{sec:compute} we report on a new
best known packing-covering lattice for $d = 6$. We thereby show in particular that 
$\gamma_6 < 1.412$, revealing the phenomenon $\gamma_6
< \gamma_5$, recently suspected by Lagarias and Pleasants
\cite{lp-2002}, Sec.~7. Note that this  
was already observed by Zong
\cite{zong-2002} (Remark 3), who showed that $\gamma(\mathsf{E}^*_6) =
\sqrt{2} < \gamma_5$.

We were not able yet to find any new best known lattices in dimensions $d
\geq 7$. The lattice $\mathsf{E}_7^*$ gives the best known
lattice in dimension~$7$. Nevertheless, because of their symmetry and
the known bounds on $\gamma_d$, Zong \cite{zong-2002} (Conjecture~3.1)
made the following conjectures: $\mathsf{E}_8$ and Leech lattice
$\Lambda_{24}$ are optimal in their dimensions. In \cite{sv-2004b} we
showed that the Leech lattice gives a locally optimal lattice
packing-covering constant. It is an open question whether or not the
root lattice $\mathsf{E}_8$ gives a locally optimal lattice
packing-covering constant as well. The corresponding values of the
lattice packing-covering constant are shown in Table~4. In
dimension~$7$ and between dimensions~$8$ and~$24$ we do not yet know
enough to state any serious conjectures.  The exact value for
the smallest known lattice packing-covering constant in dimension~$6$
is~$2\sqrt{2\sqrt{798}-56}$ (see Section~\ref{sec:compute}).

\begin{center}
\begin{tabular}{|c|c|c|}
\hline
$\mathbf d$ & \textbf{lattice} & \textbf{lattice packing-covering constant} 
$\mathbf \gamma_d$ \\ 
\hline
 $6$ & $L_6^{pc}$        & $2\sqrt{2\sqrt{798}-56}\approx 1.411081$ \\
 $7$ & $\mathsf{E}^*_7$  & $\sqrt{7/3} \approx 1.527525$ \\
 $8$ & $\mathsf{E}_8$    & $\sqrt{2} \approx 1.414213$ \\
$24$ & $\Lambda_{24}$    & $\sqrt{2} \approx 1.414213$ \\
\hline
\end{tabular}
\par\smallskip
\textbf{\textsf{Table 4.}} Lattice packing-covering records.
\end{center}


\section{Voronoi's Reduction Theory}
\label{sec:voronoireduction}

The general task of a reduction theory for PQFs is to give a
fundamental domain for $\MS^d_{>0} / \GL_d(\Z)$. This is a subset
which behaves like $\MS^d_{>0} / \GL_d(\Z)$ up to boundary
identifications. There are many different reduction theories,
connected with names like Lagrange, Gau\ss, Hermite, Korkine,
Zolotareff, Minkowski, Voronoi, and others (see \cite{so-1985}). In
this section we describe the reduction theory developed by
Voronoi in \cite{voronoi-1908}. It is based on Delone triangulations.

\subsection{Secondary Cones of Delone Triangulations}
\label{ssec:secondary}

Let $Q \in \MS^d_{>0}$ be a PQF. A polytope $L = \conv\{\vec{v}_1,
\ldots, \vec{v}_n\}$, with $\vec{v}_1, \ldots, \vec{v}_n$ in $\Z^d$,
is called a \textit{Delone polytope} of~$Q$ if there exists a $\vec{c}
\in \R^d$ and a real number $r \in \R$ with $Q[\vec{v}_i - \vec{c}] =
r^2$ for all $i = 1, \ldots, n$, and for all other lattice points
$\vec{v} \in \Z^d \backslash \{\vec{v}_1, \ldots, \vec{v}_n\}$ we have
strict inequality $Q[\vec{v} - \vec{c}] > r^2$. The set of all Delone
polytopes
\[
\Del(Q) = \{L : \mbox{$L$ is a Delone polytope of $Q$}\}
\] 
is called the \textit{Delone subdivision} (or \textit{$L$-partition})
of~$Q$. A \textit{Delone triangulation} is a Delone subdivision that
consists of simplices only. For historical reasons we use the letter
$L$ to denote Delone polytopes (see \cite{voronoi-1908} and
\cite{delone-1938}).

The Delone subdivision of a PQF is a periodic polytopal subdivision
of~$\R^d$.  We say that two Delone polytopes~$L$ and~$L'$ are
\textit{equivalent} if there is a~$\vec{v} \in \Z^d$ with $L' = L +
\vec{v}$. Given a Delone subdivision $\MD$ of $\R^d$, the set of PQFs
with Delone subdivision $\MD$ forms the \textit{secondary cone}
\[
\SC(\MD) = \{Q \in \sdo : \Del(Q) = \MD\}.
\]
In the literature it is often referred to as the \textit{$L$-type
domain} of~$Q\in\SC(\MD)$. We prefer the term ``secondary cone''
because of the close connection of Voronoi's reduction theory to the
theory of secondary polytopes which we will point out in
Section~\ref{ssec:secondarypolytopes}.

Let $Q$ be a PQF whose Delone subdivision is a triangulation
of~$\R^d$. In the following we will describe the secondary cone
of~$\Del(Q)$. For this, let $L = \conv\{\vec{v}_1, \ldots,
\vec{v}_{d+1}\}$ and $L' = \conv\{\vec{v}_{2}, \ldots,
\vec{v}_{d+2}\}$ be two $d$-dimensional Delone simplices of~$Q$
sharing the common facet $F = \conv\{\vec{v}_2, \ldots,
\vec{v}_{d+1}\}$. Let $\alpha_1, \ldots, \alpha_{d+2}$ be real numbers
with $\alpha_1 = 1$, $\sum_{i=1}^{d+2} \alpha_i = 0$ and
$\sum_{i=1}^{d+2} \alpha_i \vec{v}_i = \vec{0}$ (hence $\alpha_{d+2} >
0$). The \textit{regulator} of the pair of adjacent simplices $(L,
L')$ is the linear form $\varrho_{(L,L')}(Q') = \sum_{i=1}^{d+2}
\alpha_i Q'[\vec{v}_i]$, $Q' \in \sd$. In particular, note that the
regulator solely depends on the points $\vec{v}_1, \ldots,
\vec{v}_{d+2}$, that $\varrho_{(L,L')}(Q) > 0$, and that
$\varrho_{(L,L')} = \varrho_{(L+\vec{v},L'+\vec{v})}$ for all $\vec{v}
\in \Z^d$.  One can describe $\SC(\Del(Q))$ by linear inequalities
coming from the (finitely many) regulators of~$\Del(Q)$:

\begin{proposition} (\cite{voronoi-1908}, \S 77)
\label{prop:secondarycones}
Let $Q$ be a PQF whose Delone subdivision is a triangulation. The
secondary cone of the Delone triangulation $\Del(Q)$ is the
full-dimensional open polyhedral cone
\[
\SC(\Del(Q)) = \{Q' \in \MS^d : 
\mbox{$\varrho_{(L,L')}(Q') > 0$, 
for all pairs $(L,L')$ of adjacent simplices of $\Del(Q)$}\}.
\]
\end{proposition}

Note on the one side that $\SC(\Del(Q))$ is contained in $\sdo$ by
definition and on the other side that the linear inequalities given by
the regulators imply that a quadratic form which satisfies them is
positive definite.

\begin{example}
\label{ex:principal}
As a first example and because of its importance for the lattice
problems introduced in Section~\ref{sec:notations}, we describe the
Delone subdivision of \textit{Voronoi's principal form of the first
type} $Q[\vec{x}] = d \sum x_i^2 - \sum x_i x_j$, which is associated
to the lattice ${\mathsf A}^*_d$, in greater detail. The Delone
subdivision of~$Q$ is a triangulation and can be described as follows:
Let $\vec{e}_1, \ldots, \vec{e}_d$ be the standard basis vectors
of~$\Z^d$, and set $\vec{e}_{d+1} = -\vec{e}_1 - \cdots -
\vec{e}_d$. For a permutation $\pi \in {\mathsf S}_{d+1}$ we define
the $d$-dimensional simplex $L_{\pi}$ by
\[
L_{\pi} = \conv\{\vec{e}_{\pi(1)}, \vec{e}_{\pi(1)} +
\vec{e}_{\pi(2)}, \ldots, \vec{e}_{\pi(1)} + \cdots + \vec{e}_{\pi(d+1)}\}.
\] 
The set of simplices $\{L_{\pi} + \vec{v} : \vec{v} \in \Z^d, \pi \in
{\mathsf S}_{d+1}\}$ defines a triangulation of $\R^d$ which we from
now on denote by~$\MD_1$.  The full-dimensional cells containing the
origin are $L_\pi$, $\pi \in {\mathsf S}_{d+1}$. Two simplices
$L_{\pi}$ and $L_{\pi'}$ have a facet in common if and only if $\pi$
and $\pi'$ differ by a single transposition of two adjacent
positions. The automorphism group of~$\MD_1$ is isomorphic to the
permutation group~${\mathsf S}_{d+1}$.  The 
{\em star} of the origin is illustrated in Figure~3.  This consists of all
Delone polytopes containing the origin.

\begin{figure}
\begin{center}
\begin{center}
\fbox{ 
\unitlength1cm
\begin{picture}(5,5)
\put(2.5,0){\vector(0,1){5}}
\put(0,2.5){\vector(1,0){5}}
\put(0.5,2.5){\circle*{0.1}}
\put(0.5,0.5){\circle*{0.1}}
\put(2.5,0.5){\circle*{0.1}}
\put(2.5,2.5){\circle*{0.1}}
\put(2.5,4.5){\circle*{0.1}}
\put(4.5,2.5){\circle*{0.1}}
\put(4.5,4.5){\circle*{0.1}}
\linethickness{1pt}
\thicklines
\put(0.5,2.5){\line(0,-1){2}}
\put(0.5,2.5){\line(1,0){4}}
\put(2.5,4.5){\line(0,-1){4}}
\put(0.5,2.5){\line(1,1){2}}
\put(0.5,0.5){\line(1,1){4}}
\put(0.5,0.5){\line(1,0){2}}
\put(2.5,4.5){\line(1,0){2}}
\put(2.5,4.5){\line(0,-1){2}}
\put(2.5,0.5){\line(1,1){2}}
\put(4.5,4.5){\line(0,-1){2}}
\put(3.6,3){$L_{id}$}
\put(2.8,1.7){$L_{(23)}$}
\put(2.7,3.7){$L_{(12)}$}
\put(0.7,1.7){$L_{(13)}$}
\put(1.4,3){$L_{(123)}$}
\put(1.5,1){$L_{(132)}$}
\end{picture}
}
\end{center}
\smallskip
\textbf{\textsf{Figure 3.}} The triangulation~$\MD_1$ in dimension $d = 2$.
\end{center}
\end{figure}

The secondary cone of $\MD_1$ is (see \cite{voronoi-1908}, \S 102--104)
\[
\SC(\MD_1) = \{Q \in \MS^d : \mbox{$q_{ij} < 0$ for $i\neq j$ 
and $\sum_{i} q_{ij} > 0$ for $j=1,\dots,d$ }\}.
\] 
Its topological closure $\overline{\SC(\MD_1)}$ is called
\textit{Voronoi's principal domain of the first type}.  \hfill
\qedsymbol
\end{example}

It was shown by Voronoi (\cite{voronoi-1908}, \S 97) that the
topological closures of the secondary cones gives a facet-to-facet
tessellation of $\MS^d_{\geq 0}$.  By a theorem of Gruber and Ryshkov
\cite{gr-1989} we even have a face-to-face tessellation because
``facet-to-facet implies face-to-face''. This means that a face
sharing relative interior points with another face of its dimension
coincides with this face, whenever this property holds for the facets
(faces of co-dimension $1$).

\subsection{Bistellar Operations}
\label{ssec:bistellar}

Now, given a secondary cone, how do we find its neighbors, that is,
those secondary cones sharing a facet with the given one?

An answer can be given by taking a closer look at the so-called
\textit{repartitioning polytopes} (introduced by Ryshkov and
Baranovskii in \cite{rb-1976}), which are ``hidden'' in the definition
of the regulators. Repartitioning polytopes are $d$-dimensional
Delone polytopes having a representation as the convex hull of two
Delone simplices sharing a common facet. Thus, repartitioning
polytopes have $d + 2$ vertices.

Generally, $d$-dimensional polytopes with $d + 2$ vertices have the
special property that there are exactly two different ways to
triangulate them: Let $V$ be a set of $d + 2$ points which affinely
spans $\R^d$.  Let $\sum_{\vec{v} \in V} \alpha_{\vec{v}} \vec{v} =
\vec{0}$, $\sum_{\vec{v} \in V} \alpha_{\vec{v}} = 0$, be an affine
relation between these points. There exist exactly two triangulations
of $\conv V$: $\MT_+(V,\vec{\alpha})$ with $d$-simplices $\conv(V
\backslash \{\vec{v}\})$, $\alpha_{\vec{v}} > 0$, and
$\MT_-(V,\vec{\alpha})$ with $d$-simplices $\conv(V \backslash
\{\vec{v}\})$, $\alpha_{\vec{v}} < 0$ (see for example
\cite{gkz-1994}, Ch.~7, Prop.~1.2).

Let $\MD$ be a Delone triangulation of~$\R^d$ and let $F$ be a $(d -
1)$-dimensional cell of~$\MD$. Then, $F$ is contained in two
simplices $L$~and~$L'$ of~$\MD$. By~$V$ we denote the set of vertices
of $L$~and~$L'$, $V = \vertex L \cup \vertex L'$. By $\vec{\alpha}$ we
denote an affine relation between the points in~$V$. The
$(d-1)$-dimensional cell~$F$ is called a \textit{flippable facet} of
the triangulation~$\MD$ if one of the triangulations
$\MT_+(V,\vec{\alpha})$ or $\MT_-(V,\vec{\alpha})$ is a subcomplex
of~$\MD$. If $F$ is a flippable facet of $\MD$ and we replace the
subcomplex $\MT_+(V,\vec{\alpha})$ by $\MT_-(V,\vec{\alpha})$
[respectively $\MT_-(V,\vec{\alpha})$ by $\MT_+(V,\vec{\alpha})$],
then we get a new triangulation. This replacement is called
\textit{bistellar operation} or \textit{flip}.

Notice that non-flippable facets do exist and that performing a
bistellar operation in a Delone triangulation does not necessarily
produce a Delone triangulation. Both phenomena occur starting from
dimension~$4$.

Nevertheless, the facets of $\overline{\SC(\MD)}$ correspond exactly to those
bistellar operations of~$\MD$ which yield new Delone triangulations. A
$(d-1)$-dimensional cell $L \cap L' \in \MD$ is a flippable facet
whenever the corresponding regulator $\varrho_{(L,L')}$ gives a
facet-defining hyperplane of~$\overline{\SC(\MD)}$ (see
\cite{voronoi-1908}, \S 87--88). This is clear since the
repartitioning polytope~$\conv(L \cup L)$ is a Delone polytope of the
PQFs lying in the relative interior of the facet given
by~$\varrho_{(L,L')}$.

Let $\mathbf{F}$ be a facet of~$\overline{\SC(\MD)}$. We describe how
the Delone triangulation~$\MD$ changes if we vary a PQF continuously:
We start from the interior of~$\overline{\SC(\MD)}$, then we move
towards a relative interior point of~$\mathbf{F}$ and finally we go
infinitesimally further, leaving $\overline{\SC(\MD)}$. In every
repartitioning polytope $V = \conv(L \cup L')$ where $L, L'$ is a pair
of adjacent simplices whose regulator defines~$\mathbf{F}$, i.e.\ the
linear span of $F$ satisfies $\lin \mathbf{F} = \{Q \in \MS^d :
\varrho_{(L,L')}(Q) = 0\}$, we perform a bistellar operation. This
gives a new triangulation~$\MD'$, which is again a Delone
triangulation. The two secondary cones $\overline{\SC(\MD)}$ and
$\overline{\SC(\MD')}$ have the complete facet~$\mathbf{F}$ in
common. We say that $\MD$ and $\MD'$ are \textit{bistellar
neighbors}. In \cite{voronoi-1908} \S 91--96, Voronoi computes the
secondary cone of~$\MD'$ explicitly and shows that
$\overline{\SC(\MD')}$ has dimension $\frac{d(d+1)}{2}$.

\subsection{Main Theorem of Voronoi's Reduction Theory}
\label{ssec:maintheorem}

By constructing bistellar neighbors we could produce infinitely many
Delone triangulations starting from the Delone triangulation~$\MD_1$
of Voronoi's principal form of the first type (a part of the infinite
\textit{flip graph} of two-dimensional Delone triangulations is given
in Figure~4). Many of these will not be essentially new, because the
group $\GL_d(\Z)$ is acting on the set of Delone subdivisions by
\hbox{$(A, \MD) \mapsto A\MD$} and it is acting on the set of
secondary cones by $(A, \SC) \mapsto A^t \SC A$ for $A \in \GL_d(\Z)$.
We are only interested in the orbits of these group actions and there
are only finitely many, as shown by Voronoi \cite{voronoi-1908} \S 98
(see also Deza, Grishukhin and Laurent in \cite{dl-1997}, Chapter
13.3). Altogether we get:

\begin{theorem} (\textit{Main Theorem of Voronoi's Reduction Theory})
\label{th:mainvoronoi}
\finishline The secondary cone of a Delone triangulation is a
full-dimensional, open polyhedral cone in $\MS^d_{\geq 0}$. The
topological closures of secondary cones of Delone triangulations give
a face-to-face tiling of $\MS^d_{\geq 0}$.  Two secondary cones share
a facet if and only if they are bistellar neighbors. The group
$\GL_d(\Z)$ acts on the tiling, and under this group action there are
only finitely many inequivalent secondary cones.
\end{theorem}

\begin{figure}
\begin{center}
\begin{center}
\fbox{
\includegraphics[width=12cm]{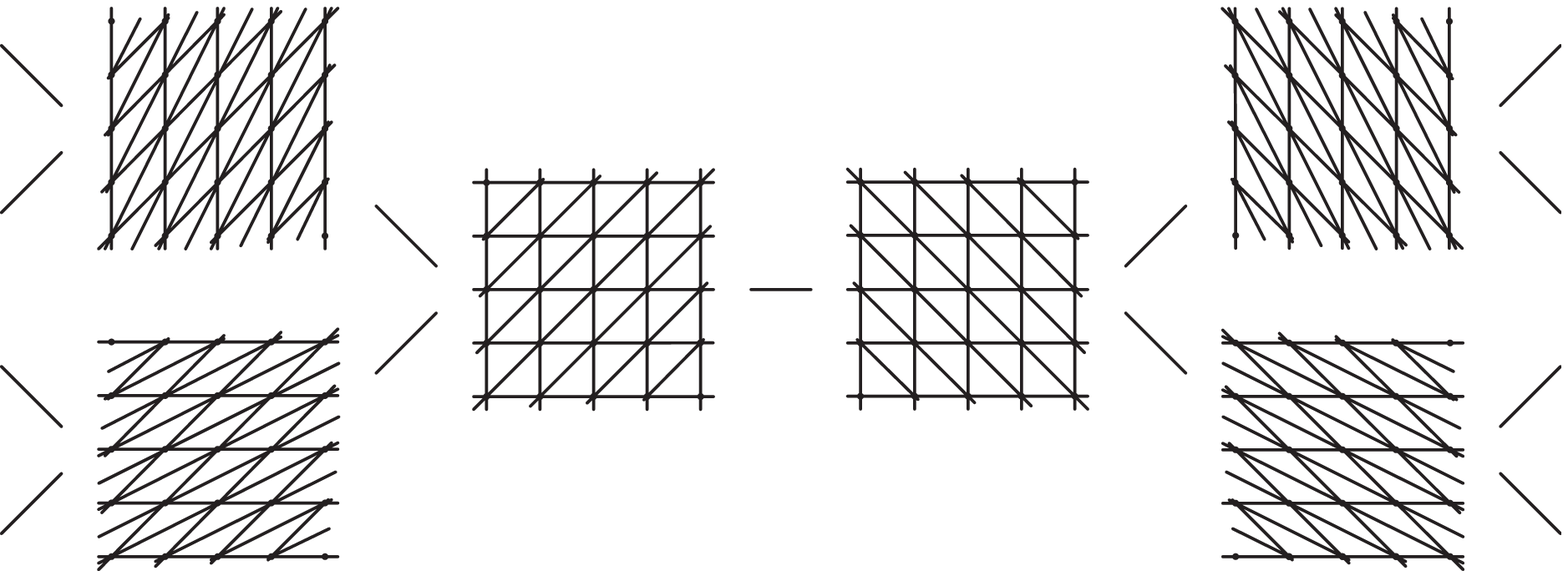}
}
\end{center}
\smallskip
\textbf{\textsf{Figure 4.}} The graph of two-dimensional Delone triangulations.
\end{center}
\end{figure}

The main theorem translates into Algorithm 1 which enumerates all
inequivalent Delone triangulations in a given dimension. We developed
the program \texttt{scc} (\textit{secondary cone cruiser}) 
which is an implementation of Algorithm~1. 
The interested reader can download \texttt{scc}  from
our web page \cite{sv-2005}. Using our
implementation we succeeded in reproducing the classification of all inequivalent
Delone triangulations up to dimension~$5$. Table~5 shows the numbers.

\begin{figure}
\begin{center}
\begin{center}
\fbox{
\begin{minipage}{14cm}
\begin{flushleft}
\smallskip
\textbf{Input:} Dimension $d$.\\
\textbf{Output:} Set~$\MR$ of all inequivalent $d$-dimensional Delone triangulations.\\
\smallskip
$T \leftarrow \{\MD_1\}$, where $\MD_1$ is the Delone triangulation described in Example~\ref{ex:principal}.\\
$\MR \leftarrow \emptyset$.\\
\textbf{while} there is a $\MD \in T$ \textbf{do}\\
\hspace{2ex} $T \leftarrow T \backslash \{\MD\}$. $\MR \leftarrow \MR \cup \{\MD\}$.\\
\hspace{2ex} compute the regulators of~$\MD$.\\
\hspace{2ex} compute the facets $F_1, \ldots, F_n$ of $\overline{\SC(\MD)}$.\\
\hspace{2ex} \textbf{for} $i = 1, \ldots, n$ \textbf{do}\\
\hspace{2ex}\hspace{2ex} compute the bistellar neighbor $\MD_i$ of $\MD$ which is defined by $F_i$.\\
\hspace{2ex}\hspace{2ex}  \textbf{if} $\MD_i$ is not equivalent to a Delone triangulation in $\MR \cup \{\MD_1, \ldots, \MD_{i-1}\}$ \textbf{then}\\
\hspace{2ex}\hspace{2ex}\hspace{2ex} $T \leftarrow T \cup \{\MD_i\}$.\\
\hspace{2ex}\hspace{2ex} \textbf{end if}\\
\hspace{2ex} \textbf{end for}\\
\textbf{end while}\\
\end{flushleft}
\end{minipage}
}
\end{center} 
\smallskip
\textbf{\textsf{Algorithm 1.}} Enumeration of all inequivalent Delone triangulations.
\end{center}
\end{figure}

We are not the first to compute this classification. Voronoi
performed the classification of all inequivalent Delone triangulation
up to dimension~$4$ in his memoir \cite{voronoi-1908}. In
\cite{br-1973} and \cite{rb-1976} Ryshkov and Baranovskii reported on
$221$ inequivalent Delone triangulation in dimension~$5$. However,
they missed one type which was found by Engel
\cite{engel-1998}. Grishukhin and Engel \cite{eg-2002} undertook the
non-trivial task of identifying the Delone triangulation missing in the
list of Ryshkov and Baranovskii. There they also report on several
errors in both lists. Our computations confirm that the number 
of Delone triangulations in dimension~$5$ is $222$. Beginning with dimension~$6$
the number of inequivalent Delone triangulations starts to explode.
We found more than $250,000$ inequivalent triangulations before our
implementation stopped because of memory reasons. We do not know how
many inequivalent triangulations we have to expect in dimension~$6$
but we are quite certain that we only saw a small fraction of them. No
non-trivial bounds on the number of inequivalent Delone triangulations
for general dimension~$d$ are known.

\begin{center}
\begin{tabular}{|l|c|c|c|c|c|}
\hline
\textbf{Dimension} & $1$ & $2$ & $3$ & $4$ & $5$\\
\hline
\textbf{\# Delone triangulations} & $1$ & $1$ & $1$ & $3$ & $222$ \\
\hline
\end{tabular}
\par\smallskip
\textbf{\textsf{Table 5.}} Numbers of inequivalent Delone triangulations.
\end{center}

\subsection{Degeneracy}
\label{ssec:degenerations}

Until now we have only dealt with Delone \textit{triangulations} of
positive \textit{definite} quadratic forms. Let us look at possible
degenerations --- Delone subdivisions of positive semidefinite
quadratic forms --- and find out how they fit into the theory
developed so far.

Let $Q$ be a positive semidefinite quadratic form which is
arithmetically equivalent to $\left(\begin{smallmatrix}0 & 0\\ 0 &
Q'\end{smallmatrix}\right)$ where $Q'$ is positive definite. Then, we
can use the definition of Delone subdivision almost literally; we only
have to replace ``polytope'' by ``polyhedron'' (a polyhedron is the
intersection of finitely many half spaces; a polytope is a bounded
polyhedron).

Delone subdivisions are limiting cases of triangulations. Their
secondary cones occur on the boundaries of full-dimensional secondary
cones of Delone triangulations. Let $\MD$~and~$\MD'$ be two Delone
subdivisions. We say $\MD$ is a \textit{refinement} of $\MD'$ if every
Delone polytope of~$\MD$ is a subset of some Delone polytope
of~$\MD'$.  The following proposition, which seems to be folklore, shows
that the relation between refinements, secondary cones and sums of
positive semidefinite quadratic forms is very natural. One can find a
proof for example in \cite{vallentin-2003} (Prop.~2.6.1).

\begin{proposition}
\label{prop:delonedegeneration}
Let~$\MD$ be a Delone triangulation.
\begin{enumerate}
\item A positive semidefinite quadratic form $Q$ lies
in~$\overline{\SC(\MD)}$ if and only if~$\MD$ is a refinement of~$\Del(Q)$.
\item If two positive semidefinite quadratic forms $Q$ and $Q'$ both
lie in~$\overline{\SC(\MD)}$, then $\Del(Q+Q')$ is a common refinement
of~$\Del(Q)$ and $\Del(Q')$.
\end{enumerate}
\end{proposition}

Therefore, the classification of all inequivalent Delone subdivisions
is equivalent to the classification of all inequivalent secondary
cones. This has been done up to dimension~$5$. The $1$- and
$2$-dimensional cases are trivial; the $3$-dimensional case goes back
to Federov, who classified all polytopes which tile $3$-dimensional
space by translates in 1885. Delone \cite{delone-1929} (later
corrected by Stogrin \cite{stogrin-1973}) found $51$ of the $52$
Delone subdivisions in dimension~$4$.  Recently, Engel
\cite{engel-2000} reported that there are $179,372$ inequivalent
five-dimensional Delone subdivisions. It is possible to verify his
result by enumerating the faces of the $222$ secondary cones of Delone
triangulations in dimension $5$ up to equivalence. Since we do not
need Engel's result for our application we did not verify it. Again,
the number in dimension~$6$ is not known and will be much
larger. Table~6 summarizes this discussion.

\begin{center}
\begin{tabular}{|l|c|c|c|c|c|}
\hline
\textbf{Dimension} & $1$ & $2$ & $3$ & $4$ & $5$\\
\hline
\textbf{\# Delone subdivisions} & $1$ & $2$ & $5$ & $52$ & $179,372$\\
\hline
\end{tabular}
\par\smallskip
\textbf{\textsf{Table 6.}} Numbers of inequivalent Delone subdivisions.
\end{center}

\subsection{Generalized Secondary Polytopes}
\label{ssec:secondarypolytopes}

Triangulations of discrete point sets have attracted many researchers
in recent years. They have many applications, for example in computational
geometry, optimization, algebraic geometry, topology, etc. One main
tool to understand the structural behavior of triangulations of finite
point sets is the theory of secondary polytopes invented by Gel'fand,
Kapranov and Zelevinsky (\cite{gkz-1994}, Ch.~7). 

Let $\MA = \{\vec{a}_1, \ldots, \vec{a}_n\} \subseteq \R^d$ be a
finite set of points. Let $w : \MA \to \R$ be a map that assigns to
every point in~$\MA$ a weight. The set of weight maps forms a vector
space over~$\R$ which we denote by~$\R^{\MA}$. A lifting map~$l : \MA
\to \R^d \times \R$, $l(\vec{a}_i) = (\vec{a}_i, w(\vec{a}_i))$ is
defined by $w$ which lifts each point $\vec{a}_i \in \MA$ by its
weight~$w(\vec{a}_i)$. A subdivision of the convex polytope~$\conv
\MA$ is induced by~$l$: We take the convex hull of the lifted
points~$\conv l(\MA)$ and project its lower faces as seen from
$(\vec{0}, -\infty)$ back down onto~$\R^d$. A subdivision that can be
obtained in this manner is called a \textit{regular
subdivision}. Delone subdivisions (or more precisely Delone
subdivisions of finitely many points) are regular subdivisions since
the underlying positive semidefinite quadratic form is the weight
function. This view on Delone subdivisions was introduced by Brown
\cite{brown-1979} and by Edelsbrunner and Seidel \cite{es-1986}.

Let~$\MT$ be a regular triangulation of~$\conv \MA$. We may ask what
are the weight functions which define~$\MT$. What is the
\textit{secondary cone} of~$\MT$ in the parameter space~$\R^{\MA}$?
As in Voronoi's reduction theory it turns out that the secondary
cone of~$\MT$ is a full-dimensional open polyhedral cone.  The
topological closures of the secondary cones of all regular
triangulations tile the space~$\R^{\MA}$ face-to-face. The tiling
is called \textit{secondary fan} of~$\MA$. If two secondary cones have
a facet in common, then the corresponding regular triangulations
differ by a bistellar operation in exactly one ``repartitioning
polytope'' (in this context it is a polytope with $d+2$ vertices
without the condition of being a Delone polytope) that is defined by
the facet. The faces in the secondary fan~$\MA$ are in a one-to-one
correspondence to regular subdivisions in essentially the same way
we discussed in Section~\ref{ssec:degenerations} for Delone
subdivisions.

So far we have seen that the theory of regular subdivisions of finite
point sets and the theory of Delone subdivisions of the lattice~$\Z^d$
can be developed analogously, but there are also differences: 
The parameter space for regular 
subdivisions is the vector space~$\R^{\MA}$, 
while for Delone subdivisions it is the
pointed cone~$\MS^d_{\geq 0}$. Groups play an important role for
Delone subdivisions. The group $\Z^d$ is acting on Delone subdivisions
by translations. On the set of secondary cones the group~$\GL_d(\Z)$
is acting.

If we order all regular subdivisions of~$\conv \MA$ by refinement we
get a poset. This poset has a very nice combinatorially structure as
proved by Gel'fand, Kapranov and Zelevinsky : There exists a polytope
--- the \textit{secondary polytope} $\Sigma(\MA)$ of $\MA$ --- whose
normal fan equals the secondary fan of~$\MA$. So the refinement poset
is anti-isomorphic to the face lattice of the secondary
polytope. Regular triangulations are in one-to-one correspondence to
the vertices, two regular triangulations differ by a bistellar
operation if and only if their vertices are connected by an edge,
etc. 

Recently, a similar combinatorial structure lurking behind the
refinement poset of $\Z^d$-periodic subdivisions (this poset contains
the poset of Delone subdivisions) has been described by Alexeev
\cite{alexeev-2002}, Sec.~5.11 and~5.12 . He gives an
\textit{unbalanced} and a \textit{balanced} version of these
\textit{generalized secondary polytopes}, where the latter one is
invariant with respect to the group action of $\GL_d(\Z)$.


\section{Convex Optimization with LMI Constraints}
\label{sec:convopt}

In this section we introduce determinant maximization problems, which
are convex programming problems with linear matrix inequality
constraints. In a sense they are equivalent to the better known
semidefinite programming problems. For both classes efficient
algorithms and implementations are available. In
Section~\ref{sec:algorithms} we will see how Problem~\ref{problem:cov}
and Problem~\ref{problem:packcov} can naturally be formulated as a
finite number of determinant maximization problems.

Following Vandenberghe, Boyd, and Wu \cite{vbw-1998} we say that a
\textit{determinant maximization problem} is an optimization problem
of the form
\begin{equation}
\label{eqn:primal-problem}
\boxed{
\begin{array}{ll}
\mbox{\textbf{minimize}} & \vec{c}^t\vec{x}-\log\det G(\vec{x})\\
\mbox{\textbf{subject to}} & \mbox{$G(\vec{x})$ is a positive definite matrix},\\
& \mbox{$F(\vec{x})$ is a positive semidefinite matrix}.
\end{array}
}
\end{equation}
The optimization vector is $\vec{x} \in \R^d$, the objective vector
is $\vec{c} \in \R^d$, and $G : \R^d \to \R^{m \times m}$ and $F :
\R^d \to \R^{n \times n}$ are affine maps:
\[
\begin{array}{rcl}
G(\vec{x}) & = & G_0 + x_1 G_1 + \cdots + x_d G_d,\\
F(\vec{x}) & = & F_0 + x_1 F_1 + \cdots + x_d F_d,
\end{array}
\] 
where $G_i \in \R^{m \times m}$ and $F_i \in 
\R^{n \times n}$, for $i = 0,\ldots, d$, are symmetric matrices.  
In the sequel we will write
$G(\vec{x}) \succ 0$ and $F(\vec{x}) \succeq 0$ for the \textit{linear
matrix inequalities} defining the constraints of the determinant
maximization problem. As a special case, our
formulation reduces to a semidefinite programming problem whenever
$G(\vec{x})$ is the identity matrix for all $\vec{x} \in \R^d$.

Currently there exist two different types of algorithms --- ellipsoid
and interior-point methods --- which efficiently solve semidefinite
programming problems. They can approximate the solution of a
semidefinite programming problem within any specified accuracy and run
in polynomial time if the instances are ``well-behaved''. For more
information on the exciting topic of semidefinite programming the
interested reader is referred to the vast amount of literature which
to a great extend is available on the World Wide Web. A good starting
point is the 
web page\footnote{\href{http://www-user.tu-chemnitz.de/~helmberg/semidef.html}
{\texttt{http://www-user.tu-chemnitz.de/\~{}helmberg/semidef.html}}}
of Christoph Helmberg.

Nesterov and Nemirovskii \cite{nn-1994} developed a framework for the
design of efficient interior-point algorithms for general and
specific classes of convex programming problems. There (\S 6.4.3),
they also showed that the determinant maximization problem can be
transformed into a semidefinite programming problem by a
transformation which can be computed in polynomial time. Nevertheless
it is faster to solve the determinant maximization problem directly.
Vandenberghe, Boyd, Wu \cite{vbw-1998} and independently Toh
\cite{toh-1999} gave interior-point algorithms for the determinant
maximization problem. Both algorithms fit into the general framework
of Nesterov and Nemirovskii.  For our implementation we use the
software package
\texttt{MAXDET}\footnote{\href{http://www.stanford.edu/~boyd/MAXDET.html}{\texttt{http://www.stanford.edu/\~{}boyd/MAXDET.html}}}
of Wu, Vandenberghe, and Boyd as a subroutine.

One nice feature of determinant maximization problems is that there is
a duality theory similar to the one of linear programming, which
allows one to compute \textit{certificates} for a range in which the optimum 
is attained, the so called \textit{duality gap}.  The dual problem of problem
\eqref{eqn:primal-problem} is (see \cite{vbw-1998}, Sec.~3)

\begin{equation}  \label{eqn:dual-problem}
\boxed{
\begin{array}{ll}
\mbox{\textbf{maximize}} & \log\det W - \trace(G_0 W) - \trace(F_0 Z) + m\\
\mbox{\textbf{subject to}} & \mbox{$\trace(G_i W) + \trace(F_i Z) = \vec{c}_i$ for $i = 1, \dots, d$,}\\
& \mbox{$W \succ 0$, $Z \succeq 0$.}
\end{array}
}
\end{equation}
Here, $W\in \R^{m\times m}$ and $Z\in\R^{n\times n}$ are symmetric matrices
of the same size as the $G_i$ and $F_i$ respectively.

The knowledge of a vector $\vec{x}$ with $G(\vec{x}) \succ 0$ and
$F(\vec{x}) \succeq 0$, and of a pair $(W, Z)$, with $\trace(G_i
W)+\trace(F_i Z)=\vec{c}_i$ for $i = 1, \dots, d$, and $W \succ 0$, and
$Z \succeq 0$, gives upper and lower bounds for the optimal value
$p^*$ of \eqref{eqn:primal-problem} by
\[
\log\det W - \trace(G_0 W) - \trace(F_0 Z) + m
\leq p^* \leq
\vec{c}^t\vec{x} - \log\det G(\vec{x}).
\]
If $\vec{x}$, $W$ and $Z$ have rational entries only, we can compute
lower and upper bounds of $p^*$ that are mathematical rigorous.  For
the proofs in Section~\ref{sec:compute} we have therefore developed a
program \texttt{rmd} (\textit{rigorous MAXDET}) which calls
\texttt{MAXDET} to find floating point values of a vector $\vec{x}$
and a pair $(W,Z)$ satisfying the constraints. Then it approximates
$\vec{x}$ and $(W,Z)$ by rational numbers and checks whether these
rational approximations satisfy the constraints. In such a case they
guarantee a certified duality gap. The interested reader can download
\texttt{rmd} from our web page \cite{sv-2005}.


\section{An LMI Constraint for the Inhomogeneous Minimum}
\label{sec:lmi}

In this section we will give a linear matrix inequality in the
parameters $(q_{ij})$ of a PQF $Q = (q_{ij})$, which is satisfied if
and only if the inhomogeneous minimum of~$Q$ is bounded by a constant,
say~$\mu(Q) \leq 1$.  For this it is crucial to observe that
$\vec{x}^t Q \vec{y}$ is a linear expression in the parameters
$(q_{ij})$.  The PQF $Q = (q_{ij})$ defines the inner product of a
Euclidean space $(\R^d, (\cdot, \cdot))$ by $(\vec{x}, \vec{y}) =
\vec{x}^tQ\vec{y}$.

From the article \cite{ddrs-1970}, \S3, of Delone, Dolbilin, Ryshkov
and Stogrin we can extract the following proposition which will be
central in our further discussion.

\begin{proposition}
\label{prop:lmi}
Let $L = \conv\{\vec{0}, \vec{v}_1, \ldots, \vec{v}_d\}\subseteq
\R^d$ be a $d$-dimensional simplex. Then the radius of the circumsphere 
of~$L$ is at most~$1$ with respect to $(\cdot, \cdot)$ 
if and only if the following linear matrix
inequality (in the parameters $q_{ij}$) is satisfied:
\[
\BR_L(Q) =
\begin{pmatrix}
4 & (\vec{v}_1, \vec{v}_1) & (\vec{v}_2, \vec{v}_2) & \ldots &
(\vec{v}_d, \vec{v}_d)\\
(\vec{v}_1, \vec{v}_1) & (\vec{v}_1, \vec{v}_1) & (\vec{v}_1, \vec{v}_2) & \ldots &
(\vec{v}_1, \vec{v}_d)\\
(\vec{v}_2, \vec{v}_2) & (\vec{v}_2, \vec{v}_1) & (\vec{v}_2,
\vec{v}_2) & \ldots & (\vec{v}_2, \vec{v}_d)\\
\vdots & \vdots & \vdots & \ddots & \vdots\\
(\vec{v}_d, \vec{v}_d) & (\vec{v}_d, \vec{v}_1) & (\vec{v}_d, \vec{v}_2) & \ldots &
(\vec{v}_d, \vec{v}_d)\\
\end{pmatrix}
\succeq 0.
\]
\end{proposition}

In \cite{ddrs-1970} Delone et al.\ used this proposition as the key
ingredient for showing that the set of PQFs which determine a
circumsphere of~$L$ with radius of at most~$R$ is convex and bounded
(see Proposition~\ref{prop:properties}).  Because of the importance of
the proposition and for the convenience of the reader we give a short
proof here.

\begin{proof}
We make use of Cayley-Menger determinants. The Cayley-Menger
determinant of $d+1$ points $\vec{v}_0, \ldots, \vec{v}_d$ with
pairwise distances $\dist(\vec{v}_i,\vec{v}_j)$ is
\[
\CM(\vec{v}_0, \ldots, \vec{v}_d) =
\begin{vmatrix}
0 & 1 & \ldots & 1 \\ 1 & \dist(\vec{v}_0,\vec{v}_0)^2 & \ldots &
\dist(\vec{v}_0,\vec{v}_d)^2\\
\vdots & \vdots                       & \ddots & \vdots                       \\
1      & \dist(\vec{v}_d,\vec{v}_0)^2 & \ldots & \dist(\vec{v}_d,\vec{v}_d)^2\\
\end{vmatrix}.
\] 
The squared circumradius $R^2$ of the simplex~$L$ equals (see for example
\cite{berger-1987} Prop.~9.7.3.7)
\[
R^2 = -\frac{1}{2}
\cdot
\frac
{
\det
\left(
\dist(\vec{v}_i, \vec{v}_j)^2
\right)_{0 \leq i,j \leq d}
}
{\CM(\vec{v}_0, \ldots, \vec{v}_d)}.
\] 
Replacing $\dist(\vec{x}, \vec{y})^2$ by $(\vec{x}, \vec{x}) -
2(\vec{x}, \vec{y}) + (\vec{y}, \vec{y})$, using $\vec{v}_0 =
\vec{0}$, and performing elementary transformations of the
determinants turns the above formula into
\begin{equation}
\label{eq:circumradius}
R^2 = 
-\frac{1}{4} 
\cdot
\frac
{
\begin{vmatrix}
0 & (\vec{v}_1, \vec{v}_1) & (\vec{v}_2, \vec{v}_2) & \ldots &
(\vec{v}_d, \vec{v}_d)\\
(\vec{v}_1, \vec{v}_1) & (\vec{v}_1, \vec{v}_1) & (\vec{v}_1, \vec{v}_2) & \ldots &
(\vec{v}_1, \vec{v}_d)\\
\vdots & \vdots & \vdots & \ddots & \vdots\\
(\vec{v}_d, \vec{v}_d) & (\vec{v}_d, \vec{v}_1) & (\vec{v}_d, \vec{v}_2) & \ldots &
(\vec{v}_d, \vec{v}_d)\\
\end{vmatrix}
}
{\det\left((\vec{v}_i, \vec{v}_j)\right)_{1 \leq i,j \leq d}}.
\end{equation}
The inequality $R \leq 1$ is equivalent to
\[
4 \cdot {\det\left((\vec{v}_i, \vec{v}_j)\right)_{1 \leq i,j \leq d}} 
+ 
\begin{vmatrix}
0 & (\vec{v}_1, \vec{v}_1) & (\vec{v}_2, \vec{v}_2) & \ldots &
(\vec{v}_d, \vec{v}_d)\\
(\vec{v}_1, \vec{v}_1) & (\vec{v}_1, \vec{v}_1) & (\vec{v}_1, \vec{v}_2) & \ldots &
(\vec{v}_1, \vec{v}_d)\\
\vdots & \vdots & \vdots & \ddots & \vdots\\
(\vec{v}_d, \vec{v}_d) & (\vec{v}_d, \vec{v}_1) & (\vec{v}_d, \vec{v}_2) & \ldots &
(\vec{v}_d, \vec{v}_d)\\
\end{vmatrix}
\geq 0,
\]
which is equivalent to
\[
\begin{vmatrix}
4 & (\vec{v}_1, \vec{v}_1) & (\vec{v}_2, \vec{v}_2) & \ldots &
(\vec{v}_d, \vec{v}_d)\\
(\vec{v}_1, \vec{v}_1) & (\vec{v}_1, \vec{v}_1) & (\vec{v}_1, \vec{v}_2) & \ldots &
(\vec{v}_1, \vec{v}_d)\\
\vdots & \vdots & \vdots & \ddots & \vdots\\
(\vec{v}_d, \vec{v}_d) & (\vec{v}_d, \vec{v}_1) & (\vec{v}_d, \vec{v}_2) & \ldots &
(\vec{v}_d, \vec{v}_d)\\
\end{vmatrix}
\geq 0.
\] 
The minors in the lower right are all determinants of Gram matrices
and therefore non-negative.  Hence, the matrix $\BR_L(Q)$ is positive
semidefinite.
\end{proof}

\begin{corollary} \label{cor:brdet}
For any $d$-dimensional simplex $L \subseteq \R^d$ with vertex at
$\vec{0}$ and $Q\in\sdo$ we have
\[
\left|\BR_L(Q)\right|\geq 0 
\quad
\Longleftrightarrow
\quad
\BR_L(Q)
\succeq 0
.
\]
\end{corollary}

\begin{example}
Let us compute the matrix linear inequality $\BR_L(Q) \succeq 0$, $Q =
\left(\begin{smallmatrix}q_{11} & q_{21}\\q_{21} & q_{22}
\end{smallmatrix}\right)$, for the two-dimensional simplex $L =
\conv\{\binom{0}{0}, \binom{1}{0}, \binom{1}{1}\}$. We have
\begin{eqnarray*}
\BR_L(Q) & = &
\begin{pmatrix}
4                         & q_{11}          & q_{11} + 2q_{21} + q_{22} \\
q_{11}                    & q_{11}          & q_{11} + q_{21}           \\
q_{11} + 2q_{21} + q_{22} & q_{11} + q_{21} & q_{11} + 2q_{21} + q_{22} \\
\end{pmatrix}
\\
& = &
\begin{pmatrix}
4 & 0 & 0\\
0 & 0 & 0\\
0 & 0 & 0\\
\end{pmatrix}
+ q_{11}
\begin{pmatrix}
0 & 1 & 1\\
1 & 1 & 1\\
1 & 1 & 1\\
\end{pmatrix}
+ q_{21}
\begin{pmatrix}
0 & 0 & 2\\
0 & 0 & 1\\
2 & 1 & 2\\
\end{pmatrix}
+ q_{22}
\begin{pmatrix}
0 & 0 & 1\\
0 & 0 & 0\\
1 & 0 & 1\\
\end{pmatrix}.
\end{eqnarray*}
\hfill \qedsymbol
\end{example}

\noindent For a Delone polytope $L$ other than a simplex the
circumradius is less than or equal to $1$ if and only if it is less
than or equal to $1$ for some $d$-dimensional simplex $L'$ with
vertices in $\vertex L$.  Therefore we set $\BR_{L}(Q)=\BR_{L'}(Q)$ in
this case.  Since a block matrix is semidefinite if and only if the
blocks are semidefinite, we have the following proposition which
allows us to express the constraint ``$\mu(Q)\leq 1$'' locally by a
single LMI.

\begin{proposition}  
\label{propos:mu}
Let $Q = (q_{ij}) \in \sdo$ be a PQF and let $L_1, \ldots, L_n$ be a
representative system of all inequivalent $d$-dimensional Delone
polytopes in $\Del(Q)$. Then
\[
\mu(Q)\leq 1
\quad\Longleftrightarrow\quad
\begin{pmatrix}
\boxed{\BR_{L_1}(Q)} & 0                   & 0 & \ldots & 0\\
0                   & \boxed{\BR_{L_2}(Q)} & 0 & \ldots & 0\\
\hdotsfor{5}\\
0                   & 0                   & 0 & \ldots & \boxed{\BR_{L_n}(Q)}
\end{pmatrix}
\succeq 0.
\]
\end{proposition}


\section{Algorithms}
\label{sec:algorithms}

In this section we present algorithms which in theory solve the
lattice covering problem and the lattice packing-covering problem in
any dimension~$d$.

Our algorithm for the lattice covering problem computes all locally
optimal lattice coverings of a given dimension. These are only
finitely many because we will see that for every fixed Delone
triangulation $\MD$ there exists at most one PQF lying in the
topological closure of the secondary cone of~$\MD$ giving a locally
optimal covering density. So, we fix a Delone triangulation and try to
find the PQF which minimizes the density function in the
topological closure of the secondary cone of the fixed Delone
triangulation. We will formulate this restricted lattice covering
problem as a determinant maximization problem.

Our algorithm for the lattice packing-covering problem operates
similarly. For every Delone triangulation we have to solve a
semidefinite programming problem.

\subsection{Solving the Lattice Covering Problem}
\label{ssec:solvinglcp}

Recall that the covering density of a PQF~$Q$ in $d$~variables is
given by $\Theta(Q) = \sqrt{\frac{\mu(Q)^d}{\det Q}} \cdot
\kappa_d$. Scaling of $Q$ by a positive real number~$\alpha$ leaves
$\Theta$ invariant. Consequently, we can restrict our attention to
those PQFs $Q$ with $\mu(Q) = 1$. Hence, we solve the lattice
covering problem if we solve the optimization problem
\[
\boxed{
\begin{array}{ll}
\mbox{\textbf{maximize}} & \det(Q)\\
\mbox{\textbf{subject to}} & \mbox{$Q$ is a positive definite matrix},\\
& \mbox{$\mu(Q) = 1$},\\
\end{array}
}
\] 
where the optimization variables $q_{ij}$ are the entries of the
PQF~$Q$. The major disadvantage of this optimization problem is that
the second constraint is not expressible as a convex condition in the
variables $q_{ij}$ and that the problem has many local maxima. A
locally optimal solution is also called a \textit{locally optimal
lattice covering}.

We will circumvent this by splitting the original problem into a
finite number of determinant maximization problems. For every Delone
triangulation~$\MD$ we solve the optimization problem
\[
\boxed{
\begin{array}{ll}
\mbox{\textbf{maximize}} & \det(Q)\\
\mbox{\textbf{subject to}} & \mbox{$Q \in \overline{\SC(\MD)}$},\\
& \mbox{$\mu(Q) \leq 1$}.\\
\end{array}
}
\] 
The relaxation of no longer requiring $\mu(Q) = 1$ in the third
constraint does not give more optimal solutions because if $Q$ 
satisfies the constraints, then so does $\frac{1}{\mu(Q)}Q$. Now, we have to show
that this is indeed a determinant maximization problem. We have seen
in Proposition~\ref{prop:secondarycones} that the first constraint can
be expressed with inequalities linear in~$q_{ij}$. The constraint
$\mu(Q) \leq 1$ can be expressed by a linear matrix inequality as in
Proposition~\ref{propos:mu}.

The optimization vector $\vec{x} \in \R^{d(d+1)/2}$ is the vector of
coefficients of $Q$. The linear matrix inequality $G(\vec{x}) \succ 0$
is given by $G(\vec{x}) = Q$. We encode the two other constraints $Q
\in \overline{\SC(Q)}$ and $\mu(Q) \leq 1$ by block matrices in the
linear matrix inequality $F(\vec{x}) \succeq 0$. For any linear
inequality which is needed to describe the secondary cone we have a $1
\times 1$ block matrix.  For any inequivalent $d$-dimensional
simplex $L \in \MD$ we have the $(d+1)\times(d+1)$ block matrix
$\BR_L(Q)$.

\subsection{Solving the Lattice Packing-Covering Problem}
\label{ssec:solvingslpc}

Along the same lines as above we formulate the lattice
packing-covering problem as a finite number of semidefinite
programming problems.

Recall that the packing-covering constant of a PQF~$Q$ is
$\gamma(Q) = 2 \cdot\sqrt{\frac{\mu(Q)}{\lambda(Q)}}$.  Since
$\gamma$ is homogeneous we can again assume $\mu(Q) =
1$, and the lattice packing-covering problem is
equivalent to the following optimization problem.

\[
\boxed{
\begin{array}{ll}
\mbox{\textbf{maximize}} & \lambda(Q)\\
\mbox{\textbf{subject to}} & \mbox{$Q$ is a positive definite matrix},\\
& \mbox{$\mu(Q) = 1$}.\\
\end{array}
}
\] 
A locally optimal solution of the optimization problem is called a
\textit{locally optimal lattice packing-covering}.

From the previous discussion we know how to deal with the constraint
$\mu(Q) = 1$.  Now, how do we maximize $\lambda(Q)$? We say that
$\vec{v} \in \Z^d \backslash \{\vec{0}\}$ is a \textit{shortest
vector} of $Q$ if $Q[\vec{v}] = \lambda(Q)$.  A theorem of Voronoi
(see \cite{voronoi-1908} \S 55, or \cite{cs-1988} Ch.~21, Th.~10)
implies that a shortest vector $\vec{v}$ gives the edge $[\vec{0},
\vec{v}]$ in the Delone subdivision of~$Q$. In a fixed Delone
subdivision~$\MD$ there are only finitely many (at most $2(2^d - 1)$,
see \cite{voronoi-1908} \S 55) edges of the form $[\vec{0},
\vec{v}]$. We can maximize $\lambda(Q)$ for $Q \in
\overline{\SC(\MD)}$ as follows: We introduce a new variable $C$ which
we want to maximize subject to the constraints $Q[\vec{v}] \geq C$
where $\vec{v}$ runs through all edges of the form $[\vec{0},
\vec{v}]$ in~$\MD$. This assures $\lambda(Q) = C$ when $C$ attains a
maximum. All these expressions are linear in the coefficients $q_{ij}$
of $Q$. Hence we have to solve the following semidefinite programming
problem for every inequivalent Delone 
triangulation~$\MD$ in order to
solve the lattice packing-covering problem:

\[
\boxed{
\begin{array}{ll}
\mbox{\textbf{maximize}} & C\\
\mbox{\textbf{subject to}} & \mbox{$Q \in \overline{\SC(\MD)}$},\\
& \mbox{$\mu(Q) \leq 1$},\\
& \mbox{$Q[\vec{v}] \geq C$, where $[\vec{0}, \vec{v}]$ is an edge in $\MD$}.\\
\end{array}
}
\]


\section{Local Optima}
\label{sec:classical}

Since we are dealing with convex optimization problems, we can extract
some structural information about uniqueness and invariance properties
of locally optimal solutions.  On the one hand these help us to
identify exact coordinates of optimal solutions. On the other hand
they allow us to decide whether we have found an isolated locally
optimal solution or not.

\subsection{Properties of Local Optima}

Most of the desired properties follow from the special structure of
sets with PQFs $Q$ attaining certain fixed values $\mu(Q)$,
$\lambda(Q)$ or $\det(Q)$. Let $L = \conv\{\vec{0}, \vec{v}_1, \ldots,
\vec{v}_d\} \subseteq \R^d$ denote a $d$-dimensional simplex and $\MD$
a Delone triangulation. Then we consider the sets
\begin{eqnarray*}
\mathbf{V}_L           & = & \{ Q \in \sdo : \left|\BR_L(Q)\right| \geq 0\},\\
\widehat{\mathbf{V}}_L & = & \{ Q \in \sdo : \left|\BR_L(Q)\right| = 0\}, \\
\mathbf{W}_\MD         & = & \bigcap_{L \in \MD} \mathbf{V}_L \cap \SC(\MD).
\end{eqnarray*}
Hence, by Corollary~\ref{cor:brdet} and Proposition~\ref{propos:mu}
the PQFs $Q\in \mathbf{W}_\MD$ are those with $\mu(Q)\leq 1$ in
$\SC(\MD)$. Further, for $D > 0$ and $\lambda > 0$ we consider
\begin{eqnarray*}
\mathbf{D}_D        & = & \{ Q \in \sdo : \det(Q) \geq D \},\\
\mathbf{M}_\lambda  & = & \{ Q \in \sdo : \lambda(Q) \geq \lambda \}.
\end{eqnarray*}

The following is known (see \cite{ddrs-1970}, \cite{ryshkov-1974}) and
was partially discussed in previous sections:

\begin{proposition}
\label{prop:properties}
Let $L = \conv\{\vec{0}, \vec{v}_1, \ldots, \vec{v}_d\} \subseteq \R^d$ 
be a $d$-dimensional simplex, $\MD$ a Delone triangulation,
$D > 0$ and $\lambda > 0$. Then
\begin{enumerate}
\item $\mathbf{V}_L$, and therefore $\mathbf{W}_\MD$, is convex and bounded.
\item $\widehat{\mathbf{V}}_L$ is a smooth (regular) hypersurface.
\item $\mathbf{D}_D$ is unbounded, strictly convex and has a smooth boundary.
\item $\mathbf{M}_\lambda$ is unbounded and convex and has a piecewise
linear boundary.
\item $\mathbf{D}_D$ and $\mathbf{M}_\lambda$ are invariant with
respect to the action of $\gldz$ on $\sdo$.
\end{enumerate}
\end{proposition}

As a consequence of these properties, it is not hard to derive the
following proposition (see \cite{ddrs-1970}, \cite{ryshkov-1974}).

\begin{proposition} 
\label{prop:invariance}
Let $\MD$ be a Delone triangulation.
\begin{enumerate}
\item \label{item:local-opt:1} The set of PQFs in
$\overline{\SC(\MD)}$ attaining $\min\limits_{Q\in\overline{\SC(\MD)}}
\Theta(Q)$ is a single PQF, together with all of its positive
multiples.  It is invariant with respect to the action of $\MD$'s
automorphism group
\[
\Aut(\MD) = \{U \in \GL_d(\Z) : \MD U = \MD\}.
\]
\item The set of PQFs in $\overline{\SC(\MD)}$ attaining
$\min\limits_{Q\in\overline{\SC(\MD)}} \gamma(Q)$ is convex and contains a
subset which is invariant with respect to $\Aut(\MD)$.
\end{enumerate}
\end{proposition}

A first, but non-geometric proof of property \ref{item:local-opt:1}
was given by Barnes and Dickson \cite{bd-1967}.  There they also
observed, as Ryshkov \cite{ryshkov-1974} did for the lattice
packing-covering problem, that it suffices to optimize among all PQFs
whose automorphism group contains the group~$\Aut(\MD)$.

\begin{corollary}
If $Q$ is a local optimum of the lattice covering or lattice packing-covering
problem among all PQFs in $\SC(\MD)$ whose automorphism group contains 
$\Aut(\MD)$, then $Q$ is a locally optimal solution.
\end{corollary}

Note that the statement above holds only if $Q$ lies in the
\textit{interior} of the secondary cone. For local optima on the
boundary of secondary cones we may apply the following trivial
proposition.

\begin{proposition}
\label{prop:localcriterion}
A PQF $Q$ is a locally optimal solution with respect to $\Theta$ or
$\gamma$ if and only if it is an optimal solution for all Delone
triangulations $\MD$ with $Q\in\overline{\SC(\MD)}$.
\end{proposition}

There probably exist two locally optimal solutions with respect to
$\gamma$ in dimension~$5$ which both lie on the boundary of some
secondary cones (see Section~\ref{sec:compute}, proof of
Theorem~\ref{thm:loc_pac_cov_opt}). In dimension~$5$ and lower every
locally optimal solution with respect to $\Theta$ lies in the interior
of some secondary cone. However, in higher dimension there exist
locally optimal solutions lying on the boundary. The Leech lattice
gives such an example as recently shown by the authors in
\cite{sv-2004b}.

\subsection{Rigorous Certificates}
\label{ssec:rigoruous-certificates}

Using convex optimization software to solve the covering or
packing-covering problem, we are often limited to determining a certain
certified range (see Section~\ref{sec:convopt}) in which the optimum
value is attained. To use Proposition~\ref{prop:localcriterion} it is
desirable to decide computationally whether or not the optimum is
attained on the boundary or even on a specific facet of
$\overline{\SC(\MD)}$. In many cases this is possible by using the
following proposition which is a simple consequence of the convexity
of determinant maximization problems.

\begin{proposition}
\label{prop:rigorouscriterion}
Consider the determinant maximization problem 
\begin{equation} \label{eqn:maxdet}
\begin{array}{ll}
\mbox{\textbf{minimize}} & f(\vec{x})\\
\mbox{\textbf{subject to}} & \mbox{$G(\vec{x}) \succ 0$ 
                             and $F(\vec{x}) \succeq 0$},\\
\end{array}
\end{equation}
where $f(\vec{x}) = \vec{c}^t\vec{x} - \log\det G(\vec{x})$.
\begin{enumerate}
\item 
Let $\vec{y}_1$ be the optimal solution of (\ref{eqn:maxdet}), and let
$\vec{y}_2$ be the optimal solution of (\ref{eqn:maxdet}) with the
additional constraint $\tilde{F}(\vec{x}) \succeq 0$. If $f(\vec{y}_1)
< f(\vec{y}_2)$, then $\tilde{F}(\vec{y}_2) = 0$.
\item Let $\vec{y}_1$ be an optimal solution of (\ref{eqn:maxdet})
with the additional linear constraint $\vec{n}^t \vec{x} \geq 0$, and
let $\vec{y}_2$ be an optimal solution of (\ref{eqn:maxdet}) with the
additional linear constraint $\vec{n}^t \vec{x} \leq 0$.
\begin{enumerate}
\item If $f(\vec{y}_1) < f(\vec{y}_2)$, then $\vec{n}^t \vec{y}_2 = 0$.
\item If $f(\vec{y}_1) > f(\vec{y}_2)$, then $\vec{n}^t \vec{y}_2 < 0$.
\end{enumerate}
\end{enumerate}
\end{proposition}

Let $Q_{\opt}$ denote a PQF with $\mu(Q_{\opt}) = 1$ attaining an
optimum in $\overline{\SC(\MD)}$ for $\Theta$. The lattice
packing-covering case is similar. Let $l$ and $u$ be a lower and an
upper bound for a minimum of $f(Q) = -\log\det(Q)$. Here and in the
sequel we use $\mu_{\MD}(Q)\leq 1$ as an abbreviation for the LMI in
Proposition~\ref{propos:mu} guaranteeing that the circumradius (with
respect to $Q$) of all the simplices of the triangulation $\MD$ is at
most $1$.

The first part of Proposition~\ref{prop:rigorouscriterion} gives a
sufficient criterion for $Q_{\opt}$ lying on the boundary of
$\overline{\SC(\MD)}$. Suppose the upper bound of $f$ on $\{ Q\in \sdo
: \mu_{\MD}(Q)\leq 1 \}$ is smaller than $l$. Then we have a
certificate for $Q_{\opt}\in \bd\overline{\SC(\MD)}$.

The first half of the second part of
Proposition~\ref{prop:rigorouscriterion} gives a sufficient criterion
for $Q_{\opt}$ lying on a specific facet $F$ of $\overline{\SC(\MD)}$.
Let $H$ denote a hyperplane containing $F$ and let $H^-$, $H^+$ be the
closed halfspaces containing $\overline{\SC(\MD)}$, respectively not
containing $\overline{\SC(\MD)}$.  If the upper bound of $f$ on $\{
Q\in \sdo : \mbox{$Q\in H^+$ and $\mu_{\MD}(Q)\leq 1$}\}$ is smaller
than $l$, then we have a certificate for $Q_{\opt}\in F$.

The second half of the second part of
Proposition~\ref{prop:rigorouscriterion} gives a sufficient criterion
for $Q_{\opt}$ not lying on a specific facet $F$ of
$\overline{\SC(\MD)}$.  If the lower bound of $f$ on $\{ Q\in \sdo :
\mbox{$Q \in H^+$ and $\mu_{\MD}(Q)\leq 1$}\}$ is larger than $u$,
then we have a certificate for $Q_{\opt} \not \in F$.  Clearly, such a
certificate for all facets of $\overline{\SC(\MD)}$ gives a
certificate for $Q_{\opt} \in \SC(\MD)$.

Note that these certificates give mathematical rigorous proofs when we
proceed as described at the end of Section~\ref{sec:convopt}.

\subsection{Necessary and Sufficient Conditions}

Next we assume a PQF $Q$ is given and we want to decide
computationally if it is a locally optimal solution to the lattice
covering or lattice packing-covering problem.  In
Section~\ref{sssec:interior} we consider the case when $Q$ lies in the
interior of the secondary cone of some Delone triangulation. In this
case we give necessary and sufficient conditions for $Q$ being a
locally optimal solution. The case when $Q$ lies on the boundary of
the secondary cone of some Delone triangulation is more subtle and we
deal with it in Section~\ref{sssec:boundary}.

\subsubsection{Normal Cones}
\label{sssec:normalcones}

Before considering locally optimal solution we have to do some local
analysis. Here, we shall compute the normal cones of the sets
$\mathbf{W}_\MD$, $\mathbf{M}_\lambda$ and $\mathbf{D}_D$ at a given PQF $Q$.
As a general reference to the basic concepts in convex and differential
geometry used in the sequel we refer to the book \cite{schneider-1993}
by Schneider.

We consider $\sd$ as Euclidean space with inner product $\langle \cdot
, \cdot \rangle$.  Let $H = \{ S \in \sd : \langle N, S \rangle =
\alpha \}$ be a hyperplane with normal vector $N \in
\sd\setminus\{\vec{0}\}$ and $C$ a convex set with boundary point $Q$.
Then $H$ is said to be a supporting hyperplane of $C$ at $Q$ with
outer normal vector $N$, if $C \subseteq H^-=\{ S\in\sd : \langle N, S
\rangle \leq \langle N, Q \rangle \}$.  The normal cone of $C$ at $Q$
is then given by all outer normal vectors of supporting hyperplanes at
$Q$ together with the zero vector. Clearly the normal cones are convex.  
A hyperplane $H$ is called a separating
hyperplane for two convex sets with a common boundary point $Q$, if it
is a supporting hyperplane of both sets, but with opposite outer
normal vectors.  Such a hyperplane exists if and only if the
corresponding normal cones $\MN_1$ and $\MN_2$ at $Q$ satisfy
$-\MN_1\cap \MN_2 \neq \emptyset$.

Let us now compute the normal cones of the sets $\mathbf{W}_\MD$,
$\mathbf{M}_\lambda$ and $\mathbf{D}_D$ at some PQF $Q$.

\begin{proposition}
Let $\MD$ be a Delone triangulation and $Q \in \mathbf{W}_\MD$ a PQF
with $\mu(Q) = 1$.  Then the normal cone of $\mathbf{W}_\MD$ at $Q$ is
equal to $-\cone\{g_L(Q) : L \in \MD \mbox{ and } Q \in
\widehat{\mathbf{V}}_L \}$, where $g_L(Q)=\gradient |\BR_L| (Q)$
denotes the gradient of the regular surface $\widehat{\mathbf{V}}_L$
at $Q$.
\end{proposition}

\begin{proof}
First, $-g_L(Q)\neq\vec{0}$ exists and is the unique outer normal
vector of $\mathbf{V}_L$ at $Q$, since $\widehat{\mathbf{V}}_L$ is a
regular surface defined by the polynomial equation
$\left|\BR_L(Q)\right| = 0$.  Second, we see that $\mathbf{W}_\MD$, in
a sufficiently small neighborhood of $Q$, is equal to the intersection
$\bigcap_{Q \in \widehat{\mathbf{V}}_L} \mathbf{V}_L$.  
Therefore, the normal cone of $\mathbf{W}_\MD$ at $Q$ is equal to $\cone\{-g_L(Q) :
Q\in\widehat{\mathbf{V}}_L \}$ (see \cite{schneider-1993}, Th.~2.2.1) and
the assertion follows.
\end{proof}

By Proposition \ref{prop:properties} the boundary of $\mathbf{D}_{\det(Q)}$
is smooth at $Q$, hence up to scaling there is a uniquely determined outer normal.
The following is well known (cf. Proposition \ref{prop:optimize_equidiscriminant}
and its detailed proof in \cite{vallentin-2003}, Prop.~8.2.2).

\begin{proposition} \label{prop:normal-cone1}
\label{prop:discriminantnormal}
The normal cone of $\mathbf{D}_{\det(Q)}$ at $Q$ is given by
$-\cone\{Q^{-1}\}$.
\end{proposition}

The normal cone of $\mathbf{M}_\lambda$ at a boundary point $Q$ is
determined by the shortest vectors of $Q$, that is, by those $\vec{v}
\in \Z^d$ for which the homogeneous minimum $\lambda(Q)$ is
attained. Recall that $Q[\vec{v}] = \lambda$ for a fixed $v \in \Z^d
\setminus \{\vec{0}\}$ is linear in the parameters $q_{ij}$, hence
this condition defines a hyperplane in $\sd$.

\begin{proposition} \label{prop:normal-cone2}
For $\vec{v} \in \Z^d \setminus \{\vec{0}\}$ let $V_{\vec{v}} \in
\sd$ be the normal vector of the hyperplane $\{Q' \in\sd : Q'[\vec{v}] =
\lambda(Q)\}$ with $\langle Q, V_{\vec{v}} \rangle = \lambda(Q)$. Then,
the normal cone of $\mathbf{M}_\lambda$ at $Q$ is given by $-\cone \{
V_{\vec{v}} : Q[\vec{v}] = \lambda(Q) \}$.
\end{proposition}

\begin{proof}
This follows immediately from \cite{schneider-1993}, Th.~2.2.1.
\end{proof}

Note, that a normal vector $V_{\vec{v}} $ depends on the chosen inner
product. For example for practical computations it is convenient to
identify $\sd$ with $\R^{d(d+1)/2}$ and to use the standard inner
product $\langle Q, V \rangle = \sum_{i\leq j} q_{ij} v_{ij}$ and
$V_{\vec{v}}=(v_{ij})$ with $v_{ij} = (2 - \delta_{ij}) \vec{v}_i
\vec{v}_j$.

\subsubsection{Interior Cases}
\label{sssec:interior}

Because of the convexity of $\mathbf{D}_{\det(Q)}$,
$\mathbf{M}_{\lambda(Q)}$ and $\mathbf{W}_\MD$ a separating hyperplane
at a PQF $Q$ of $\mathbf{W}_\MD$ with $\mu(Q) = 1$ and one of the
other two sets yields a necessary and sufficient condition for $Q$ to
be a locally optimal solution to either the lattice covering or the
lattice packing-covering problem.  Therefore a necessary and
sufficient condition can be derived from the normal cones at $Q$ with
respect to these sets.

\begin{proposition} 
\label{prop:localoptcrit}
Let $\MD$ be a Delone triangulation, and let $Q \in \SC(\MD)$ be a PQF
with $\mu(Q) = 1$.
\begin{enumerate}
\item
(\cite{bd-1967}) Then, $Q$ is a locally optimal solution to the lattice
covering problem if and only if
\[
\cone\{Q^{-1}\} \cap  -\cone\{g_L(Q) : \mbox{$L \in \MD$ and $
Q\in\widehat{\mathbf{V}}_L$} \} \neq \emptyset.
\]
\item Then, $Q$ is a locally optimal solution to the lattice
packing-covering problem if and only if
\[
\cone\left\{ V_{\vec{v}} : \langle Q, V_{\vec{v}} \rangle = 
\lambda(Q) \right\} 
\cap 
-\cone\{g_L(Q) : \mbox{$L \in \MD$ and $
Q\in\widehat{\mathbf{V}}_L$} \} \neq \emptyset.
\]
\end{enumerate}
\end{proposition}

\subsubsection{Boundary Cases}
\label{sssec:boundary}

Note that the foregoing proposition does not give criteria for PQFs
$Q$ on the boundary of some secondary cones. In such a case we need to
replace $-\cone\{g_L(Q) : \mbox{$L \in \MD$ and $Q \in
\widehat{\mathbf{V}}_L$}\}$ by a generalized expression. That is done
by considering for each Delone triangulation $\MD$ with $Q \in
\overline{\SC(\MD)}$ the normal cone $\MC_\MD(Q)$ of
$\overline{\mathbf{W}_\MD}$ at $Q$.  We have $\MC_\MD(Q) =
-\cone(\MN_1 \cup \MN_2)$ with outer normal vectors $\MN_1 = \cone\{g_L(Q) :
\mbox{$L \in \MD$ and $Q \in \widehat{\mathbf{V}}_L$}\}$, as considered
before, and outer normal vectors
\[
\MN_2 = \{N \in \sd : \langle N, \cdot \rangle =
\varrho_{(L,L')}(\cdot) \mbox{ for a pair $(L,L')$ of adjacent
simplices of $\MD$ with $\varrho_{(L,L')}(Q) = 0$} \}
\]
from facets of $\overline{\SC(\MD)}$ containing $Q$. Since
$\mathbf{D}_D$ has a smooth boundary a necessary condition for a local
covering optimum at $Q$ is the convexity of
$\bigcup_{Q\in\overline{\SC(\MD)}}\overline{\mathbf{W}_\MD}$ at $Q$
which is equivalent to $\bigcap_{Q\in\overline{\SC(\MD)}}\MC_\MD(Q)
\neq \emptyset$.  In case of the lattice packing-covering problem,
this convexity condition is not necessary for a local optimum.  The
following proposition however gives sufficient conditions for a local
optimum in case the convexity condition is satisfied.

\begin{proposition} 
\label{prop:local-opt-boundary}
Let $Q\in \sdo$ be a PQF with $\mu(Q) = 1$.
\begin{enumerate}
\item
Then, $Q$ is a locally optimal solution to the lattice covering problem if
and only if
\[
\cone\{Q^{-1}\} \cap \bigcap_{Q \in \overline{\SC(\MD)}} \MC_\MD(Q)
\neq \emptyset
.
\]
\item
Then, $Q$ is a locally optimal solution to the lattice packing-covering
problem if
\[
\cone\{ V_{\vec{v}} : \langle Q, V_{\vec{v}} \rangle = \lambda(Q)\} 
\cap
\bigcap_{Q\in\overline{\SC(\MD)}}\MC_\MD(Q)
\neq \emptyset
.
\]
\end{enumerate}
\end{proposition}

\begin{proof}
The set $\bigcap_{Q\in\overline{\SC(\MD)}}\MC_\MD(Q)$ contains 
all common outer normal vectors of the sets $\overline{\mathbf{W}_\MD}$
with $Q\in\overline{\SC(\MD)}$ at $Q$. In particular, 
$\bigcap_{Q\in\overline{\SC(\MD)}}\MC_\MD(Q)$ is empty if and only if
$\bigcup_{Q\in\overline{\SC(\MD)}}\overline{\mathbf{W}_\MD}$ is not
convex at $Q$.
If convex, then a separating hyperplane of
$\bigcup_{Q\in\overline{\SC(\MD)}}\overline{\mathbf{W}_\MD}$
and $\mathbf{D}_{\det(Q)}$, respectively $\mathbf{M}_{\lambda(Q)}$
at $Q$ is sufficient for a local optimum at $Q$.
The assertion then follows by Propositions \ref{prop:normal-cone1} and 
\ref{prop:normal-cone2}
\end{proof}

\subsection{Conditions for Isolated Local Optima}

Locally optimal solutions to the lattice packing-covering problem may
not be isolated optima, in contrast to solutions to the lattice
covering problem.  To determine computationally if a given locally
optimal solution $Q$ is isolated, we have to check if there exists a
segment $[Q,Q']$ which lies in the boundary of
$\mathbf{M}_{\lambda(Q)}$ and in the boundary of $\mathbf{W}_\MD$ at
$Q$.

\begin{proposition}
Let $\MD$ be a Delone triangulation and let $Q\in\SC(\MD)$ be a
locally optimal solution to the lattice packing-covering problem with
$\mu(Q) = 1$.  Then $Q$ is not an isolated local optimum if and only
if there exists an $S\in \sd$ with
\begin{enumerate}
\item $\langle V_{\vec{v}} , S \rangle \geq 0$ 
for all $\vec{v} \in \Z^d$ with $Q[\vec{v}] = 
\lambda(Q)$,
\item $\langle g_L(Q) , S \rangle \geq 0$ for all $L\in \MD$ with 
$Q\in\widehat{\mathbf{V}}_L$,
\item $h_L(Q) S = \vec{0}$ for all $L\in \MD$ with 
$Q \in \widehat{\mathbf{V}}_L$ and $\langle g_L(Q) , S \rangle = 0$.
\end{enumerate} 
\end{proposition}

\begin{proof}
The first condition says that there exists an $\epsilon_1 > 0$ so that
the segment $[Q,Q+\epsilon_1 S]$ lies in $\mathbf{M}_{\lambda(Q)}$
because $\mathbf{M}_{\lambda(Q)}$ has a piecewise linear boundary.
The first condition together with the second condition plus the local
optimality of~$Q$ are equivalent to the fact that $S$ is in the
tangent space of $\mathbf{V}_L$ at $Q$ and in the tangent space of
$\mathbf{M}_{\lambda(Q)}$ at $Q$.  Thus $Q + \epsilon_2 S$ is in the
boundary of $\mathbf{W}_\MD$ for sufficiently small $\epsilon_2 > 0$
if and only if the corresponding smooth hypersurfaces
$\widehat{\mathbf{\mathbf{V}}}_L$ have curvature $0$ in direction
$S$. This is equivalent to $S$ being an eigenvector of eigenvalue $0$
of the Hessian $h_L(Q) = \hess|\BR_L|(Q)$ which is the third
condition.  Hence, all three conditions are fulfilled if and only if
there exists a segment $[Q,Q']$ in the boundary of the two sets
$\mathbf{M}_{\lambda(Q)}$ and $\mathbf{W}_\MD$, where $Q' = Q +
\min(\epsilon_1, \epsilon_2) S$.
\end{proof}

Note again, that the condition applies only if $Q$ is in the interior
of a secondary cone. A similar sufficient condition for PQFs on the
boundary of secondary cones can be attained by replacing $g_L(Q)$ in 2
and 3 with the rays of the normal cone
$\bigcap_{Q\in\overline{\SC(\MD)}}\MC_\MD(Q)$ of
$\bigcup_{Q\in\overline{\SC(\MD)}} \overline{\mathbf{W}_\MD}$ at $Q$
(see Proposition~\ref{prop:local-opt-boundary}), if such exist.

Finally, we propose two (as far as we know) still unanswered
questions.

\begin{question}
\label{prob:local-opt}
Does there exist a locally optimal solution to the lattice
packing-covering problem which is not isolated?
\end{question}

Note that a positive answer to the following question 
would imply the existence of non-isolated locally optimal solutions.

\begin{question}
Does there exist a
(locally optimal) solution to the lattice-packing covering problem,
which does not have all the symmetries of the corresponding Delone
subdivision?
\end{question}


\section{Local Lower Bounds via Moments of Inertia}
\label{sec:moments}

In this section we give simple and efficiently computable local lower
bounds for the lattice covering density and the lattice packing-covering
constant. These bounds are ``local'' in the sense that they only apply
to those PQFs lying in the topological closure of the secondary cone
of a given Delone triangulation. For their computation we only need to
know the coordinates of the simplices of the considered Delone
triangulation. They are therefore useful tools in a heuristic search
for ``good'' PQFs.  The method goes back to Ryshkov and Delone. It is
called the method of the moments of inertia because the central idea
in its proof is analogous to the Parallel Axis Theorem of Steiner in
classical mechanics. This analogy is explained in \cite{am-2003}.

Let $P \subseteq \R^d$ be a finite set of points in $d$-dimensional
Euclidean space~$(\R^d, (\cdot, \cdot))$. We interpret the points of
$P$ as masses with unit weight. The \textit{moment of inertia} of the
points about a point~$\vec{x} \in \R^d$ is defined as $I_{\vec{x}}(P)
= \sum_{\vec{v} \in P} \dist(\vec{x}, \vec{v})^2$. The
\textit{centroid} of~$P$ (center of gravity) is given by $\vec{m} =
\frac{1}{|P|}\sum_{\vec{v} \in P} \vec{v}$.  From the equations
\[
\begin{array}{lcl}
\dist(\vec{x}, \vec{v})^2 
& = & (\vec{x}-\vec{v}, \vec{x}-\vec{v}) \\
& = & (\vec{x}-\vec{m}, \vec{x}-\vec{m}) + (\vec{m}-\vec{v}, \vec{m}-\vec{v}) 
+ 2(\vec{x}-\vec{m}, \vec{m}-\vec{v})\\ 
& = & \dist(\vec{x}, \vec{m})^2 + \dist(\vec{m}, \vec{v})^2 
+ 2(\vec{x}-\vec{m}, \vec{m}-\vec{v}),
\end{array}
\] 
and $\sum_{\vec{v} \in P} (\vec{m}-\vec{v}) = \vec{0}$, we
derive Apollonius' formula (see \cite{berger-1987} \S 9.7.6)
\begin{equation}
\label{eq:apollonius}
I_{\vec{x}}(P) = |P| \dist(\vec{x}, \vec{m})^2 + I_{\vec{m}}(P).
\end{equation}
Hence, the moment of inertia of $P$ about the centroid~$\vec{m}$ is
minimal.

If the points of $P$ form the vertices of a $d$-dimensional simplex,
then (\ref{eq:apollonius}) gives a relationship between the radius of
the circumsphere~$R$, the center of the circumsphere~$\vec{c}$, and
the moment of inertia about the centroid~$\vec{m}$ of~$P$:
\[
R^2 = \frac{I_{\vec{c}}(P)}{d+1} = 
\dist(\vec{c},\vec{m})^2 + \frac{I_{\vec{m}}(P)}{d+1}.
\]

We can compute $I_{\vec{m}}(P)$ using only the edge lengths of the
simplex~$P$. For every $\vec{w} \in \vertex P$ we have by definition
$I_{\vec{w}}(P) = \sum_{\vec{v} \in P} \dist(\vec{w},
\vec{v})^2$. Summing up and using (\ref{eq:apollonius}) gives
\[
\sum_{\vec{w} \in P} I_{\vec{w}}(P) =
\sum_{\vec{w} \in P} ((d+1)\dist(\vec{w}, \vec{m})^2 + I_{\vec{m}}(P)) =
2(d+1) I_{\vec{m}}(P).
\] 
So, we get
\begin{equation}
\label{eq:simplex_moment}
I_{\vec{m}}(P) = \frac{1}{d+1} \sum_{\{\vec{v}, \vec{w}\} \subseteq P} 
\dist(\vec{v}, \vec{w})^2.
\end{equation}

Let $\MD$ be a Delone triangulation of~$\R^d$, let $L_1, \ldots, L_n$
be the $d$-dimensional simplices of the star of a lattice point (say
for example the origin), and let $\vec{m}_i$ be the centroid of $L_i$, $i =
1, \ldots, n$.  The arithmetical mean of the moments of inertia about
the centroids of $L_i$ with respect to a PQF~$Q$ is defined to be
\[
I_{\MD}(Q) = \frac{1}{n} \sum_{i=1}^{n} I_{\vec{m}_i}(L_i),
\] 
and is called the \textit{central moment of inertia} of~$\MD$ with
respect to~$Q$. Note that we are now dealing with the inner product
given by~$Q$ with $\dist(\vec{x}, \vec{y})^2 = Q[\vec{x}-\vec{y}]$ and
that $I_{\MD}(Q)$ is linear in the parameters $q_{ij}$ of $Q$.

\begin{proposition}
\label{prop:radius_bound}
The central moment of inertia of~$\MD$ with respect to~$Q$ yields a
lower bound for the inhomogeneous minimum of~$Q$ if $\MD$ is a
refinement of $\Del(Q)$. In this case we have
\[ 
\mu(Q) \geq
\frac{1}{d+1} I_{\MD}(Q)
.
\]
\end{proposition}

\begin{proof}
Let $R_i$ be the radius and $\vec{c}_i$ be the center of the
circumsphere of the simplex $L_i$, then
\[
\begin{array}{lclcl}
\mu(Q) & = & \max\limits_{i=1,\ldots,n} R_i^2 & = &
\max\limits_{i=1,\ldots,n}\left(\dist(\vec{c}_i,\vec{m}_i)^2 +
\frac{I_{\vec{m}_i}(L_i)}{d+1}\right)\\ 
& \geq &
\max\limits_{i=1,\ldots,n} \frac{I_{\vec{m}_i}(L_i)}{d+1} & \geq &
\frac{1}{(d+1)n} \sum\limits_{i=1}^n I_{\vec{m}_i}(L_i)\\ 
& = & \frac{1}{d+1} I_{\MD}(Q).
\end{array}
\] 
\end{proof}

To find lower bounds for the lattice covering density or for the
lattice packing-covering constant, we have to minimize the linear
function $I_{\MD}$ over all PQFs with a fixed
determinant, or with a fixed homogeneous
minimum, respectively.  As in previous
sections, we consider $\sd$ as a Euclidean space with
inner product $\langle\cdot , \cdot \rangle$.

\begin{proposition}
There exists a PQF $F\in\sdo$ with $I_{\MD}=\langle F ,\cdot \rangle$.
\end{proposition}

\begin{proof}
Since $I_{\MD}$ is a linear function, there is a $F\in\sd$ with
$I_{\MD}(\cdot) = \langle F, \cdot \rangle$. For every PQF 
$Q$, we have
$\langle F, Q \rangle = I_{\MD}(Q) > 0$. Since $\MS_{>0}^d = \{Q' \in
\MS^d : \mbox{$\langle Q', Q\rangle > 0$ for all $Q \in
\MS_{>0}^d$}\}$, meaning that $\MS_{>0}^d$ is a self-dual cone, we
have $F \in \MS_{>0}^d$.
\end{proof}

Because of the preceding proposition, minimizing the linear function
$I_{\MD}$, hence $\langle F ,\cdot \rangle$, over all PQFs with a
fixed determinant $D$ is geometrically equivalent to finding the
unique PQF on the boundary of $\mathbf{D}_D$ --- the
determinant-$D$-surface --- that has a supporting hyperplane with
normal $F$ (see Proposition~\ref{prop:discriminantnormal}). Thus we
derive the following proposition (see \cite{vallentin-2003}
Prop.~8.2.2 for a detailed proof).

\begin{proposition}
\label{prop:optimize_equidiscriminant}
A linear function $f(\cdot) = \langle F, \cdot\rangle$ with $F\in\sdo$
has a unique minimum on the determinant-$D$-surface. Its value is $d
\sqrt[d]{D \det F}$ and the minimum is attained at the PQF $\sqrt[d]{D
\det F}F^{-1}$.
\end{proposition}

Now we can plug Proposition~\ref{prop:radius_bound} and
Proposition~\ref{prop:optimize_equidiscriminant} together.  This
yields a local lower bound for the covering density of a PQF.

\begin{proposition}
\label{prop:local_lower_bound}
Let $\MD$ be a Delone triangulation. Let $Q$ be a PQF for which $\MD$
is a refinement of $\Del(Q)$. Then we have a lower bound for the
covering density of~$Q$:
\[ 
\Theta(Q) \geq \Theta_*(\MD) = \sqrt{\left(\frac{d}{d+1}\right)^d \det F}
\cdot \kappa_d,
\] 
where $F$ is the positive definite matrix given by the equation
$I_{\MD}(\cdot) = \langle F, \cdot \rangle$.  Here and in the sequel
we denote the local lower bound for the Delone triangulation~$\MD$ by
$\Theta_*(\MD)$.
\end{proposition}

\begin{proof}
By Proposition~\ref{prop:optimize_equidiscriminant}, $I_{\MD}=\langle
F, \cdot \rangle$ has the unique minimum $d\sqrt[d]{D\det F}$ on the
determinant-$D$-surface.  Using this with $D = \det Q$ and applying
Proposition~\ref{prop:radius_bound} we get the lower bound
$\Theta_*(\MD)$ for~$\Theta(Q)$, because of
\[
\left(\Theta(Q)/\kappa_d\right)^2  = \frac{\mu(Q)^d}{\det Q}
\geq \left(\frac{I_{\MD}(Q)}{d+1}\right)^d / \det Q
\geq \frac{d^d \det Q \det F}{(d+1)^d \det Q}
= \left(\frac{d}{d+1}\right)^d \det F. 
\]
\end{proof}

For a local lower bound on the lattice packing-covering constant we
minimize the linear function $I_{\MD}(\cdot) = \langle F, \cdot
\rangle$ over all PQFs with a fixed homogeneous minimum. In analogy
to the above, we replace the determinant-$D$-surface (the boundary
of $\mathbf{D}_D$) by the homogeneous-minimum-$\lambda$-surface
(the boundary of $\mathbf{M}_{\lambda}$).

\begin{proposition}
\label{prop:pack-cov-lower-bound}
Let $\lambda$ be a positive number. Let $\MD$ be a Delone
triangulation and let $Q$ be a PQF for which $\MD$ is a refinement of
$\Del(Q)$.  Further, let the minimum of the linear function
$I_{\MD}(\cdot) = \langle F, \cdot \rangle$ on the
homogeneous-minimum-$\lambda$-surface be attained at a PQF
$Q_{F,\lambda}$. Then we have a lower bound for the packing-covering
constant of~$Q$:
\[ 
\gamma(Q) \geq 
2\sqrt{\frac{\langle F, Q_{F,\lambda} \rangle}{(d+1)\lambda}}.
\] 
\end{proposition}

\begin{proof}
The PQF $\frac{\lambda(Q)}{\lambda}Q_{F,\lambda}$ attains the minimum
of $\langle F, \cdot \rangle$ among PQFs in
$\textbf{M}_{\lambda(Q)}$. Applying
Proposition~\ref{prop:radius_bound} gives the desired lower bound:
\[
(\gamma(Q)/2)^2 
= \frac{\mu(Q)}{\lambda(Q)}
\geq \left(\frac{I_{\MD}(Q)}{d+1}\right) / \lambda(Q) 
\geq \frac{\langle F, \frac{\lambda(Q)}{\lambda}Q_{F,\lambda} \rangle}
{(d+1)\cdot\lambda(Q)}
= \frac{\langle F, Q_{F,\lambda} \rangle}{(d+1)\lambda}. 
\]
\end{proof}

It is difficult in practice to obtain $Q_{F,\lambda}$ and
therefore the lower bound of the proposition. Instead, if we minimize
$I_{\MD}$ on an approximation $\textbf{M}^\MD_\lambda$ of
$\textbf{M}_\lambda$, suitable for ${\MD}$, things become easier.  We
set
\[
\textbf{M}^\MD_\lambda
=
\left\{ 
Q\in\sdo : 
\mbox{$Q[\vec{v}] \geq \lambda$ for all $\vec{v} \in \Z^d$ with $[\vec{0},\vec{v}] \in \MD$} 
\right\} 
.
\]
Now finding the minimum of $I_{\MD}$ on $\textbf{M}^\MD_\lambda$
reduces to a linear program with finitely many constraints.  By
looking at equation (\ref{eq:simplex_moment}) it is easy to see
that this linear program is bounded from below since every summand
is at least $\lambda$.  Because of $\textbf{M}_\lambda \subseteq
\textbf{M}^\MD_\lambda$ we obtain the following 
practically useful lower bound for $\gamma$.

\begin{proposition}
\label{prop:pack-cov-lower-bound2}
Let $\lambda$ be a positive number. Let $\MD$ be a Delone
triangulation and let $Q$ be a PQF for which $\MD$ is a refinement of
$\Del(Q)$.  Further, let the minimum of the linear function
$I_{\MD}(\cdot) = \langle F, \cdot \rangle$ on
$\textbf{M}^\MD_\lambda$ be attained at a PQF $Q^\MD_{F,\lambda}$.
Then we have a lower bound for the packing-covering constant of~$Q$:
\[ 
\gamma(Q) \geq \gamma_*(\MD) = 
2\sqrt{\frac{\langle F, Q^\MD_{F,\lambda} \rangle}{(d+1)\cdot\lambda}}.
\] 
\end{proposition}


\section{Computational Results}
\label{sec:compute}

In Section~\ref{sec:algorithms} we developed algorithms for the
solution of the lattice covering and the lattice packing-covering
problem.  Here, we want to demonstrate that these algorithms are not
purely of theoretical interest.  We implemented the algorithms in C++,
using the package
\texttt{MAXDET}\footnote{\href{http://www.stanford.edu/~boyd/MAXDET.html}{\texttt{http://www.stanford.edu/\~{}boyd/MAXDET.html}}}
of Wu, Vandenberghe and Boyd and the package
\texttt{lrs}\footnote{\href{http://cgm.cs.mcgill.ca/~avis/C/lrs.html}{\texttt{http://cgm.cs.mcgill.ca/\~{}avis/C/lrs.html}}}
of Avis as subroutines. The interested reader can download the
implementation from our web page \cite{sv-2005}.

With the implemented algorithm we are able to determine the solution
of the lattice covering problem and the lattice packing-covering
problem in dimensions $d = 1, \ldots, 5$. This hereby reproduced the
known results. For dimension $d = 5$ this computation took on a $2$GHz
Intel Pentium computer less than $90$ minutes. By our
computations we get a conjectural list of approximations of all locally
optimal solutions. In Section~\ref{ssec:dim1to5} we give some details.

More important, with the help of a heuristic method, which we describe
in Section~\ref{ssec:heuristics}, this approach produces interesting
new lattices in dimensions $d = 6$ which give better coverings and
packing-coverings than previously known ones. We analyze them in
Section~\ref{ssec:new6} and in Section~\ref{ssec:beautification} where
we give rigorous proofs of some facts concerning these lattices which
we found experimentally. This analysis turned out to be fruitful: The
new $6$-dimensional lattices were our starting point to find new
best-known covering lattices in dimensions $7$~and~$8$. We describe
the $7$-dimensional case in Section~\ref{ssec:dim7}. The
$8$-dimensional case is explained in \cite{sv-2004b}.

\subsection{Dimensions 1, \dots, 5}
\label{ssec:dim1to5}

We described the solutions to both problems for the $1$, $2$, and
$3$-dimensional cases in Section~\ref{sec:lcp_overview}
and~\ref{sec:lspc_overview}. Since in these cases there is only one
type of Delone triangulation, these results are computationally rather
trivial.

The $4$-dimensional case requires more work because there are three
inequivalent Delone triangulations in dimension~$4$. Baranovskii
\cite{baranovskii-1965} and independently Dickson \cite{dickson-1967}
found all three inequivalent locally optimal solutions to the lattice
covering problem in dimension~$4$. Earlier Delone and Ryshkov showed
that the lattice $\mathsf{A}_4^*$ gives the best $4$-dimensional
lattice covering \textit{without} determining all locally optimal
solutions. We describe their approach in
Section~\ref{ssec:heuristics}.

For the lattice packing-covering problem the $4$-dimensional case
was resolved by Horvath \cite{horvath-1980}. He shows that there
exist, as in the covering case, three isolated, locally optimal
solutions --- one for each Delone triangulation.  It is interesting
that Horvath's optimal packing-covering lattice $\mathsf{Ho}_4$ does
not belong to the family of root lattices and their duals. An
associated PQF is 
\[
Q_{\mathsf{Ho}_4} =
\begin{pmatrix}
 2 &  1 & -1 & -1\\
 1 &  2 & -1 & -1\\
-1 & -1 &  2 &  0\\
-1 & -1 &  0 &  2
\end{pmatrix} +
\frac{1}{3}\sqrt{3}
\begin{pmatrix}
 4 &  1 & -2 & -2\\
 1 &  4 & -2 & -2\\
-2 & -2 &  4 &  1\\
-2 & -2 &  1 &  4
\end{pmatrix}.
\]
The first summand, associated to the best packing
lattice~$\mathsf{D}_4$, lies on an extreme ray of the secondary cone
belonging to $Q_{\mathsf{Ho}_4}$.

In a series of papers Ryshkov and Baranovskii solved the
$5$-dimensional lattice covering problem. In \cite{ryshkov-1973a}
Ryshkov introduced the concept of C-types. Two Delone triangulations
are of the same \textit{C-type} if their $1$-skeletons (the graph
consisting of vertices and edges of the triangulation) coincide. He
gave an algorithm to find all inequivalent C-types in any given
dimension. He computed that there are $3$ inequivalent C-types in
dimension~$4$ and that there are $76$ inequivalent C-types in
dimension~$5$.  Using this list Baranovskii and Ryshkov enumerated
$221$ (of $222$) inequivalent $5$-dimensional Delone triangulations in
\cite{br-1973}. They described the triangulations in more detail in
\cite{br-1975}. In the last paper of the series \cite{rb-1975} they
showed that the lattice~${\mathsf A}^*_5$ provides the least dense
$5$-dimensional lattice covering. In their proof they do not find all
locally optimal lattice coverings. By using estimations (the method
explained in Section~\ref{ssec:heuristics} is one of their main tools)
they merely show that all local minima exceed the covering density
of~${\mathsf A}^*_5$. As mentioned in Section~\ref{ssec:maintheorem},
Ryshkov and Baranovskii missed one Delone triangulation.  Fortunately,
this does not give a thinner lattice covering than~${\mathsf A}^*_5$.

Using our algorithm, and the techniques described in
Section~\ref{sec:classical}, to obtain certificates we produced a list
of inequivalent locally optimal lattice coverings in
dimension~$5$. This list is conjecturally complete and gives
approximations of the local optima. The computation, including all
certificates, takes about $90$ minutes on a $2$GHz Intel Pentium
computer. We can prove the following theorem rigorously.

\begin{theorem} \label{thm:loc_cov_opt}
In dimension $5$, there exist at least $216$ and at most $218$ 
inequivalent, local minima of the lattice covering density function 
$\Theta$, ranging from approximately
$2.124286$ to approximately $2.757993$. 
All of them are attained in the interior of their
secondary cones.
\end{theorem}

\begin{proof}
We take a complete list of $222$ inequivalent secondary cones,
generated by \texttt{scc} (see Section~\ref{sec:voronoireduction} and
our web page \cite{sv-2005}).  For each cone, we applied the program
\texttt{coop} (\textit{covering optimizer}), which is based on
\texttt{rmd} (\textit{rigorous MAXDET}) and available on our web page
\cite{sv-2005}.  It computes certified bounds for $\Theta$ within the
closure of a given secondary cone by solving the determinant
maximization problem described in Section~\ref{ssec:solvinglcp}.  The
bounds for the duality gap described in Section~\ref{sec:convopt} are
mathematical rigorous, since we use rational arithmetic only.  Solving
additional determinant maximization problems as described in
Section~\ref{ssec:rigoruous-certificates}, the program then tries to
obtain \textit{certificates} for the approximated local optimum to be
attained in the interior or on the boundary of the secondary cone.  As
a result, there are $216$ certified local optima attained in the
interior of secondary cones and $4$ certified non-optima attained on
the boundary.  This leaves two cases which might be local optima,
attained in the interior of their secondary cones. However, if they
are attained on the boundary (and the numerical evidence strongly
supports this) then they are not local optima, since the corresponding
secondary cones adjacent to the facets in question have a smaller
local minimum.
\end{proof} 

Note that our computations prove rigorously that the PQF of Barnes and
Trenerry \cite{bt-1972} yields the second best locally optimal lattice
covering with density of approximately $2.230117$.

\medskip

In 1986, Horvath \cite{horvath-1986} solved the lattice
packing-covering in dimension $5$. The proof is about $70$ pages long
(private communication). As in dimension $4$, Horvath's lattice
$\mathsf{Ho}_5$ is not among previously known ones. An associated PQF
is
\[
Q_{\mathsf{Ho}_5}
= 
\begin{pmatrix}
 2 &  1 &  0 & -1 & -1\\
 1 &  2 & -1 & -1 &  0\\ 
 0 & -1 &  2 &  0 &  0\\
-1 & -1 &  0 &  2 & -1\\
-1 & 0  &  0 & -1 &  2
\end{pmatrix}
+ 
\frac{1+\sqrt{13}}{6}
\begin{pmatrix}
 6 &  3 & -2 & -2 & -2\\
 3 &  6 & -2 & -3 & -2\\ 
-2 & -2 &  6 & -1 & -1\\ 
-2 & -3 & -1 &  6 &  0\\
-2 & -2 & -1 &  0 &  6
\end{pmatrix}
.
\]
Again, the first summand is associated to the best packing
lattice~$\mathsf{D}_4$. It has rank $4$ and lies on an extreme ray of
the secondary cone belonging to $Q_{\mathsf{Ho}_5}$.

By applying propositions of Section~\ref{sec:classical} we were able
to reproduce Horvath's result and moreover to attain a list of locally
optimal solutions. Again this list is conjecturally complete and gives
approximations of the local optima. Compared to the covering problem,
there are many more secondary cones that contain no locally optimal
solution.  The computation, including all certificates, also takes
about $90$ minutes on a $2$GHz Intel Pentium computer. We can prove
the following theorem rigorously.

\begin{theorem}  \label{thm:loc_pac_cov_opt}
In dimension $5$, there exist at least $47$ and at most $75$ local
minima of the lattice packing-covering constant $\gamma$, ranging from
approximately $1.449456$ to approximately $1.557564$.  At least $45$
of them are attained in the interior of their secondary cones.
\end{theorem}

\begin{proof}
As for the proof of Theorem~\ref{thm:loc_cov_opt}, we used the
complete list of $222$ inequivalent secondary cones, generated by
\texttt{scc}.  For each cone, we applied the program \texttt{pacoop}
(\textit{packing-covering optimizer}), which is also based on
\texttt{rmd} and available on our web page \cite{sv-2005}.  It
computes certified bounds for $\gamma$ within the closure of a given
secondary cone by rigorously solving the semidefinite program
described in Section~\ref{ssec:solvingslpc}.  Again the program tries
to obtain \textit{certificates} for the approximated local optimum to
be attained in the interior or on the boundary of the secondary cone.
As a result, there are $45$ certified local optima attained in the
interior of secondary cones and $147$ certified non-optima attained on
the boundary.  There are $30$ of $222$ remaining cases in which the
program \texttt{pacoop} did not give a certificate for the optimum to
be attained in the interior or the boundary of the secondary
cone. However, if all of them are attained on the boundary (and the
numerical evidence strongly supports this), there are two local optima
on the common boundary of three inequivalent secondary cones
each. Thus there exist at least $47$ local optima as claimed.
\end{proof}

\subsection{Heuristic Methods}
\label{ssec:heuristics}

Before going to dimension~$6$, let us explain one important heuristic
method, which is essential for finding the new lattices. One of the
biggest problem in finding good lattice coverings or packing-coverings
in dimension~$6$ and higher is that it is not apriori clear which
Delone triangulations admit good ones. Solving a determinant
maximization problem is a rather time consuming task and there are a
lot of inequivalent Delone triangulations to consider. So a desirable
tool is a fast computable lower bound, as described in
Section~\ref{sec:moments}.

We view the set of Delone triangulations as an undirected labeled
graph. A node represents a Delone triangulation and two nodes are
adjacent if their Delone triangulations are bistellar
neighbors. Let~$\MD$ be a Delone triangulation. We label its node by
the local lower bound~$\Theta_*(\MD)$ or~$\gamma_*(\MD)$. We can use
the labeling in two different ways.

On the one hand it is clear that if the labeling of a node is large,
then the considered Delone triangulation does not admit a good lattice
covering and lattice packing-covering respectively.  Ryshkov and
Delone \cite{dr-1963} solved the lattice covering problem in
dimension~$4$ by this method. However, the five-dimensional
lattice covering problem cannot be solved in this way. From the $222$
inequivalent Delone triangulations there are~$20$ whose local lower
bound is smaller than $\Theta({\mathsf A}_5^*)$.

On the other hand we can hope that~$\MD$ admits a good lattice
covering or packing-covering if the local lower bound is small.  In
dimension~$6$ the hope that good local lower bounds yield good lattice
coverings is partially fulfilled. We report on a typical example: We
started from the Delone triangulation of Voronoi's principal form of
the first type and took a random walk of length~$50$. Then, we found a
node labeled by~$\Theta_*(\MD)\approx 2.585149$. From this we
proceeded by taking a neighboring node having the smallest local lower
bound. By repeating this greedy strategy we resulted in a node labeled
by~$\Theta_*(\MD)\approx 2.318034$. This node is interesting for
several conjectural ``extremeness'' properties. It yields the smallest
known local lower bound and it has the largest known number of
neighbors, namely~$130$. As we will see in the next section there
exists a locally optimal lattice covering which belongs to this node
with covering density of approximately $2.466125$.  At present this is
the second best known $6$-dimensional lattice covering which is
locally optimal. Furthermore we will see that there exists a locally
optimal lattice packing-covering belonging to this node which
currently defines the best known $6$-dimensional lattice
packing-covering.

\subsection{New Six-Dimensional Lattices}
\label{ssec:new6}

We have not been able to solve the lattice covering problem in this
dimension. However we found some new interesting covering lattices. As
we reported in Section~\ref{sec:lcp_overview}, Ryshkov
\cite{ryshkov-1967} asked what is the first dimension~$d$ where
$\mathsf{A}^*_d$ does not give the least dense lattice covering.

A classification of all inequivalent Delone triangulations in
dimension~$6$ is not within reach. So we can not use our algorithm to find
the best lattice covering. Nevertheless we ran it partially to find
good lattice coverings. Since there are no other good lattice
coverings in the neighborhood of Voronoi's principal form of the first
type, we used the heuristic method described in the previous
section. In this manner we found about $100$ inequivalent secondary
cones containing lattice coverings better than ${\mathsf
A}^*_6$. Thus, we can give an answer to Ryshkov's question:

\begin{theorem} 
\label{theorem:ryshkov-answer}
Dimension $d = 6$ is the smallest dimension where the lattice
$\mathsf{A}^*_d$ does not give the least dense lattice covering.
\end{theorem} 

\begin{proof}
This theorem follows from Theorem~\ref{th:secondcovering} or from
Theorem~\ref{th:firstcovering}.
\end{proof}

Two of the new lattice coverings were strikingly good. Of course, by
the computational optimization we only got numerical estimates of
these coverings.

The now second best-known lattice covering which is locally optimal
was quite easy to find: running our heuristic method described above
finds this lattice covering in most of the trials. An approximation is
given by the PQF
\[
Q^{c2}_6 \approx 
\begin{pmatrix}
  1.9982 &  0.5270 & -0.4170 & -0.5270 &  0.5270 & -1.0541\\
  0.5270 &  1.9982 & -0.4170 & -0.5270 &  0.5270 & -1.0541\\
 -0.4170 & -0.4170 &  2.1082 & -1.0541 & -0.4170 &  0.8341\\
 -0.5270 & -0.5270 & -1.0541 &  1.9982 & -0.5270 & -0.4170\\ 
  0.5270 &  0.5270 & -0.4170 & -0.5270 &  1.9982 & -1.0541\\
 -1.0541 & -1.0541 &  0.8341 & -0.4170 & -1.0541 &  2.1082\\
\end{pmatrix}
\]
and its covering density is~$\Theta(Q^{c2}_6)\approx
2.466125$. The Delone subdivision is a triangulation and its secondary
cone has $130$~facets. The local lower bound is approximately $2.318034$.

For a while we thought that this might by the best lattice covering in
dimension~$6$, but then 
\[
Q^{c1}_6 \approx
\begin{pmatrix}
 2.0550 & -0.9424 &  1.1126 &  0.2747 & -0.9424 & -0.6153\\
-0.9424 &  1.9227 & -0.5773 & -0.7681 &  0.3651 & -0.3651\\
 1.1126 & -0.5773 &  2.0930 & -0.4934 & -0.5773 & -0.9804\\
 0.2747 & -0.7681 & -0.4934 &  1.7550 & -0.7681 &  0.7681\\
-0.9424 &  0.3651 & -0.5773 & -0.7681 &  1.9227 & -0.3651\\
-0.6153 & -0.3651 & -0.9804 &  0.7681 & -0.3651 &  1.9227\\
\end{pmatrix}
\]
with covering density $\Theta(Q^{c1}_6)\approx 2.464802$ came up.
The Delone subdivision is a triangulation and its secondary cone has
$100$ facets. The local lower bound is approximately $2.322204$. After this,
we did not find any further lattice covering records in
dimension~$6$. Furthermore, we did not find a $6$-dimensional Delone
triangulation whose secondary cone has more than $130$ facets or whose
local lower bound is less than $2.318034$.

Nevertheless, in the secondary cone of $Q^{c2}_6$ we found the PQF
\[
Q^{pc}_6 \approx
\begin{pmatrix}
 2.0088 &  0.5154 &  0.5154 & -0.5154 &  0.9778 &  0.5154\\
 0.5154 &  2.0088 &  0.5154 & -0.5154 & -0.5154 & -0.9778\\
 0.5154 &  0.5154 &  2.0088 & -0.5154 & -0.5154 &  0.5154\\
-0.5154 & -0.5154 & -0.5154 &  2.0088 & -0.9778 & -0.5154\\
 0.9778 & -0.5154 & -0.5154 & -0.9778 &  2.0088 &  0.9778\\
 0.5154 & -0.9778 &  0.5154 & -0.5154 &  0.9778 &  2.0088\\
\end{pmatrix},
\]
which gives currently the best known lattice packing-covering in
dimension $6$ with packing-covering constant $\gamma(Q^{pc}_6)\approx
1.411081$. In the next section we will examine these new PQFs in
greater detail.

\subsection{Beautification and a Unified View}
\label{ssec:beautification}

Although we found an answer to Ryshkov's question, these results are
not fully satisfying. We want to know the exact lattices and prove
rigorously that they have a good covering density and that they are
locally optimal. Even more important, we want to know an
interpretation of why these lattice coverings are good. To accomplish
this we collect some more data.

The automorphism group of $\Del(Q^{c2}_6)$ has order $3840$ and the
one of $\Del(Q^{c1}_6)$ has order $240$. With the knowledge of the
groups we were able to compute the extreme rays of both secondary
cones \cite{dv-2004}. The secondary cone $\mathbf{C}_1 =
\overline{\SC(\Del(Q^{c2}_6))}$ has $7,145,429$ and the cone
$\mathbf{C}_2 = \overline{\SC(\Del(Q^{c1}_6))}$ has $2,257,616$
extreme rays. Both contain an extreme ray associated to the lattice
$\mathsf{E}^*_6$ given for example by the PQF 
\[
Q_{\mathsf{E}_6^*} =
\begin{pmatrix}
 4 & 1 & 2 & 2 &-1 & 1\\
 1 & 4 & 2 & 2 & 2 & 1\\
 2 & 2 & 4 & 1 & 1 & 2\\
 2 & 2 & 1 & 4 & 1 & 2\\
-1 & 2 & 1 & 1 & 4 & 2\\
 1 & 1 & 2 & 2 & 2 & 4
\end{pmatrix}.
\]  
After transforming $Q^{c1}_6$ and $Q^{c2}_6$ by integral
unimodular transformations we can assume that the Delone
triangulations of the two PQFs are refinements of the Delone
subdivision $\Del(Q_{\mathsf{E}^*_6})$. This Delone subdivision was
investigated in different contexts (see \cite{worley-1987},
\cite{cs-1991}, \cite{mp-1995} and \cite{baranovskii-1991}). We
briefly review the main results:

\begin{proposition}
In the star of the origin are $720$~full-dimensional
$6$-dimensional Delone polytopes, and the automorphism group of
$Q_{\mathsf{E}^*_6}$ acts transitively on these full-dimensional
polytopes. Each polytope is the convex hull of three regular triangles
lying in three pairwise orthogonal affine planes. Each polytope has
$9$ vertices, $27$ facets, and three different triangulations, where
each triangulation consists of nine $6$-dimensional simplices. The
covering density is $\Theta(Q_{\mathsf{E}^*_6}) = \frac{8}{9\sqrt{3}}
\cdot \kappa_6 \approx 2.652071$.
\end{proposition}

\begin{proof}
Except for the possible refining triangulations, all this data is
well-known, see for example \cite{cs-1991}, \textit{Summary for}
$\textsf{E}_6^*$.

To describe the triangulations we introduce coordinates. A Delone
polytope of $Q_{\mathsf{E}_6^*}$ is similar to the polytope $P =
\conv\{\vec{u}_1, \vec{u}_{\omega}, \vec{u}_{\overline{\omega}},
\vec{v}_1, \vec{v}_{\omega}, \vec{v}_{\overline{\omega}}, \vec{w}_1,
\vec{w}_{\omega}, \vec{w}_{\overline{\omega}}, \}$ where $\vec{u}_1 =
\vec{e}_1$, $\vec{u}_{\omega} = -\frac{1}{2}\vec{e}_1 +
\frac{\sqrt{3}}{2}\vec{e}_2$, $\vec{u}_{\overline{\omega}} =
-\frac{1}{2}\vec{e}_1 - \frac{\sqrt{3}}{2}\vec{e}_2$, $\vec{v}_1 =
\vec{e}_3$, $\vec{v}_{\omega} = -\frac{1}{2}\vec{e}_3 +
\frac{\sqrt{3}}{2}\vec{e}_4$, $\vec{v}_{\overline{\omega}} =
-\frac{1}{2}\vec{e}_3 - \frac{\sqrt{3}}{2}\vec{e}_4$, $\vec{w}_1 =
\vec{e}_5$, $\vec{w}_{\omega} = -\frac{1}{2}\vec{e}_5 +
\frac{\sqrt{3}}{2}\vec{e}_6$, $\vec{w}_{\overline{\omega}} =
-\frac{1}{2}\vec{e}_5 - \frac{\sqrt{3}}{2}\vec{e}_6$. The three
possible triangulations are given by the set of nine $6$-dimensional
simplices
\begin{eqnarray*}
\MT_{\vec{u}} & = &
\left\{\conv\vertex P \setminus \{\vec{v}_z, \vec{w}_{z'}\} :
z,z' \in \{1, \omega, \overline{\omega}\right\}\},\\
\MT_{\vec{v}} & = &
\left\{\conv\vertex P \setminus \{\vec{u}_z, \vec{w}_{z'}\} :
z,z' \in \{1, \omega, \overline{\omega}\}\right\},\\
\MT_{\vec{w}} & = &
\left\{\conv\vertex P \setminus \{\vec{u}_z, \vec{v}_{z'}\} :
z,z' \in \{1, \omega, \overline{\omega}\}\right\}.
\end{eqnarray*}
To finish the proof one has to show that these sets indeed define
triangulations and that they are the only possible triangulations.
This can be done by a straightforward computation using the facts that
the minimal affine dependent subsets of $\vertex P$ are
\[
\{
\vec{v}_1, \vec{v}_{\omega}, \vec{v}_{\overline{\omega}},
\vec{w}_1, \vec{w}_{\omega}, \vec{w}_{\overline{\omega}}
\},
\{
\vec{u}_1, \vec{u}_{\omega}, \vec{u}_{\overline{\omega}},
\vec{w}_1, \vec{w}_{\omega}, \vec{w}_{\overline{\omega}}
\},
\{
\vec{u}_1, \vec{u}_{\omega}, \vec{u}_{\overline{\omega}},
\vec{v}_1, \vec{v}_{\omega}, \vec{v}_{\overline{\omega}}
\},
\]
that the facets of $P$ are $\conv \vertex P \setminus \{\vec{u}_z,
\vec{v}_{z'}, \vec{w}_{z''}\}$ where $z,z',z'' \in
\{1,\omega,\overline{\omega}\}$, and by applying the following
proposition.

\begin{proposition} (\cite{rambau-1997}, Prop.~2.2)
\label{prop:combtria}
Let $\MA \subseteq \R^d$ be a finite point set, and let $\MT$ be a set
of $d$-dimensional simplices with vertices in $\MA$. The set $\MT$
defines a triangulation of the polytope $\conv \MA$ if and only if the
following two conditions hold:
\begin{enumerate}
\item For all $S,S' \in \MT$ there exists a minimal affine dependency
$\sum_{\vec{a} \in \MA} \lambda_{\vec{a}} \vec{a} = \vec{0}$ with
$\sum_{\vec{a} \in \MA} \lambda_{\vec{a}} = 0$ so that $\{\vec{a} \in
\MA : \lambda_{\vec{a}} > 0\} \subseteq S$ and $\{\vec{a} \in \MA :
\lambda_{\vec{a}} < 0\} \subseteq S'$.
\item For all $S \in \MT$ and for every $(d-1)$-dimensional facet $F$
of $S$ there exists either a $(d-1)$-dimensional facet $F'$ of $\conv
\MA$ with $F \subseteq F'$ or there exists another simplex $S' \in \MT
\setminus \{S\}$ also having $F$ as a facet.
\end{enumerate}
\end{proposition}

\end{proof}

\noindent
Using this information we are able to prove that the PQF $Q^{c1}_6$ is
closely related to $Q_{\mathsf{E}^*_6}$.

\begin{theorem}
\label{th:e6d}
The PQF $Q^{c1}_6$ gives the least dense lattice covering among all
PQFs whose Delone subdivision is a refinement of the Delone
subdivision $\Del(Q_{\mathsf{E}^*_6})$.
\end{theorem}

\begin{proof}
Our proof is computational and uses a branch and cut method. We have
to show that all secondary cones of a Delone triangulation refining
$\Del(Q_{\mathsf{E}^*_6})$ do not contain a PQF with covering density
less than $\Theta(Q^{c1}_6)$. There are $40$ full-dimensional Delone
polytopes of $Q_{\mathsf{E}^*_6}$ which cannot be transformed into
each other by translations or by the map $\vec{x} \mapsto
-\vec{x}$. Since each of these Delone polytopes has $3$ possible
triangulations, the number of all periodic triangulations refining
$\Del(Q_{\mathsf{E}^*_6})$ is $3^{40}$. It is not apriori clear how to
distinguish between Delone and non-Delone triangulations, and it is
not possible to generate all $3^{40}$ triangulations. We choose a
backtracking approach instead.

We arrange partial triangulations in a tree. On every level one of the
$40$ Delone polytopes is triangulated so that we have $3^n$ nodes on
the $n$-th level. For every node $N$ we define the value
\[
\Theta_N = 
\max
\left\{
\det(Q) : 
\mbox{$Q \in \MS^6_{> 0}$, $\BR_L(Q) \succeq 0$ for all simplices $L$ of partial triangulation $N$}
\right\}
.
\] 
Obviously, $\Theta_N$ is a lower bound for the covering density of any
PQF whose Delone subdivision refines the partial triangulation~$N$. We
can compute a lower bound of this value by solving a determinant
maximization problem similar to the one in
Section~\ref{ssec:solvinglcp}. Note that this can be done rigorously
using rational arithmetic only if we proceed as described in
Section~\ref{sec:convopt}. If the lower bound is larger
than~$2.464802$, we can cut the tree at this node,
since the covering density of $Q^{c1}_6$ is less than $2.464801$ 
(see Theorem~\ref{th:firstcovering}). 

This algorithm visits exactly $432$ nodes of depth $40$. One of these
triangulations equals $\Del(Q^{c1}_6)$. All the others are equivalent
to $\Del(Q^{c1}_6)$ because the automorphism group of
$\Del(Q^{c1}_6)$, which has order $240$, is a subgroup of
$\Aut(Q_{\mathsf{E}^*_6})$ and the order of $\Aut(Q_{\mathsf{E}^*_6})$
equals $103680 = 432 \cdot 240$.
\end{proof}

This computational proof takes about two weeks on a $2$GHz Intel
Pentium computer. The source code \texttt{e6d.cc} is available from
our web page \cite{sv-2005} as part of the package \texttt{rmd}.

\medskip

The knowledge of the automorphism groups of $\Del(Q^{c1}_6)$ and
$\Del(Q^{c2}_6)$ also enables us to give a unified view on both
lattices.  We have
\[ 
\Aut(\Del(Q^{c1}_6)) \subseteq \Aut(\Del(Q^{c2}_6)) \subseteq
\Aut(Q_{\mathsf{E}^*_6}).
\] 
The automorphism group $\Aut(\Del(Q^{c2}_6))$ turns out to be the subgroup of
$\Aut(Q_{\mathsf{E}^*_6})$ stabilizing the minimal vectors
$\pm \vec{e}_1$, and $\Aut(\Del(Q^{c1}_6))$ is the intersection of
the two subgroups of $\Aut(Q_{\mathsf{E}^*_6})$ stabilizing
the minimal vectors $\pm \vec{e}_1$ and $\pm \vec{e}_2$ respectively.

The subspace $\mathbf{I}_1$ of all quadratic forms invariant under the
group $\Aut(\Del(Q^{c2}_6))$ is spanned by the PQFs
$Q_{\mathsf{E}^*_6}$ and $R_1$ (see below). At the same time,
$\mathbf{I}_1 \cap \mathbf{C}_1$ is a cone with extreme rays
$Q_{\mathsf{E}^*_6}, R_1$. By Proposition~\ref{prop:invariance},
$Q^{c2}_6$ has to lie in $\cone\{Q_{\mathsf{E}^*_6}, R_1\}$. The
subspace $\mathbf{I}_2$ of all quadratic forms invariant under the
group $\Aut(\Del(Q^{c1}_6))$ is four-dimensional. The cone
$\mathbf{I}_2 \cap \mathbf{C}_2$ has six extreme rays
$Q_{\mathsf{E}^*_6}, R_2, \ldots, R_6$, where
\[
R_1 = 
\begin{pmatrix}
12 & 3 & 6 & 6 & -3 & 3\\
 3 & 7 & 4 & 4 &  3 & 2\\ 
 6 & 4 & 8 & 3 &  1 & 4\\
 6 & 4 & 3 & 8 &  1 & 4\\ 
-3 & 3 & 1 & 1 &  7 & 3\\
 3 & 2 & 4 & 4 &  3 & 7
\end{pmatrix}
\quad
R_2 =
\begin{pmatrix}
0 & 0 & 0 & 0 & 0 & 0\\
0 & 5 & 2 & 2 & 3 & 1\\
0 & 2 & 4 & 0 & 2 & 2\\
0 & 2 & 0 & 4 & 2 & 2\\
0 & 3 & 2 & 2 & 5 & 3\\
0 & 1 & 2 & 2 & 3 & 5
\end{pmatrix}
\quad
R_3 =
\begin{pmatrix}
6 & 4 & 4 & 4 & 0 & 2\\
4 &11 & 6 & 6 & 5 & 3\\
4 & 6 & 8 & 3 & 3 & 4\\
4 & 6 & 3 & 8 & 3 & 4\\
0 & 5 & 3 & 3 & 7 & 4\\
2 & 3 & 4 & 4 & 4 & 7
\end{pmatrix}
\]
\[
R_4 =
\begin{pmatrix}
3 & 2 & 2 & 2 & 0 & 1\\
2 & 3 & 2 & 2 & 1 & 1\\
2 & 2 & 4 & 1 & 1 & 2\\
2 & 2 & 1 & 4 & 1 & 2\\
0 & 1 & 1 & 1 & 3 & 2\\
1 & 1 & 2 & 2 & 2 & 4
\end{pmatrix}
\quad
R_5 =
\begin{pmatrix}
7 & 3 & 4 & 4 &-1 & 2\\ 
3 &17 & 8 & 8 & 9 & 4\\
4 & 8 &12 & 3 & 5 & 6\\
4 & 8 & 3 &12 & 5 & 6\\
-1& 9 & 5 & 5 &13 & 7\\ 
2 & 4 & 6 & 6 & 7 &12
\end{pmatrix}
\quad
R_6 =
\begin{pmatrix}
9 & 6 & 6 & 6 & 0 & 3\\
6 & 9 & 6 & 6 & 3 & 3\\
6 & 6 & 8 & 4 & 2 & 4\\
6 & 6 & 4 & 8 & 2 & 4\\
0 & 3 & 2 & 2 & 5 & 3\\
3 & 3 & 4 & 4 & 3 &  6
\end{pmatrix}
\]
Note that $R_2$ lies in $\mathbf{I}_1$. Altogether, this yields
the ``picture'' in dimension~$21$ given in Figure 5. 

\begin{figure}
\begin{center}
\begin{center}
\fbox{
\unitlength1cm
\begin{picture}(6,8)
\put(0,0){\includegraphics[width=5cm]{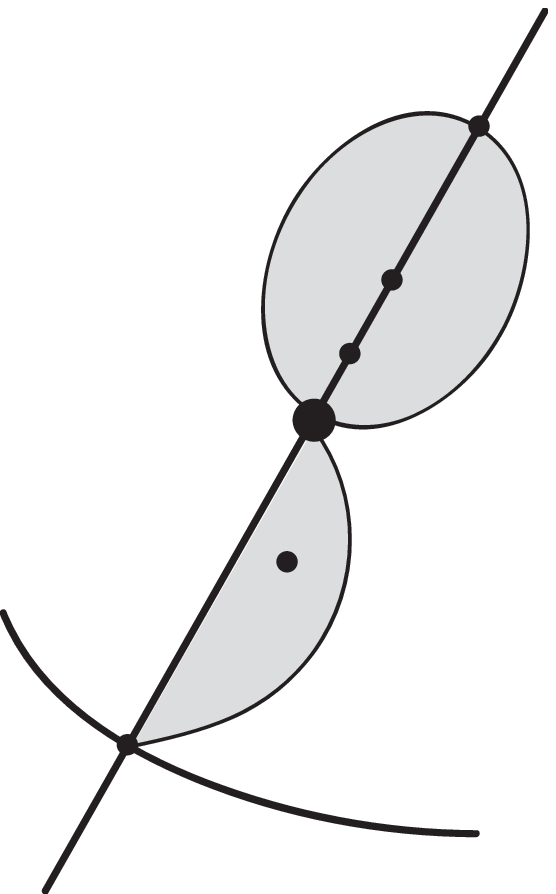}}
\put(2,4.1){$Q_{\mathsf{E}^*_6}$}
\put(4.5,6.8){$R_1$}
\put(0.8,1.6){$R_2$}
\put(2.55,3.25){$Q^{c1}_6$}
\put(3.7,5.3){$Q^{c2}_6$}
\put(2.5,4.8){$Q^{pc}_6$}
\put(3.1,6){$\mathbf{C}_1$}
\put(2.2,2.2){$\mathbf{C}_2$}
\put(2.6,7.6){invariant subspace $\mathbf{I}_1$}
\put(3.2,0.9){boundary of $\MS^6_{\geq 0}$}
\end{picture}
}
\end{center}
\smallskip
\textbf{\textsf{Figure 5.}} Unified view on $Q_{\mathsf{E}^*_6}$,  
$Q^{c1}_6$,$Q^{c2}_6$ and $Q^{pc}_6$.
\end{center}
\end{figure}

Let us finally try to find the exact coordinates of the PQFs.  This is
easy for $Q^{c2}_6$. We know that we can scale $Q^{c2}_6$ so that
$Q^{c2}_6 = Q_{\mathsf{E}^*_6} + x R_1$, for some $x \in \R_{\geq
0}$. Now, the exact finding of $x$ boils down to finding roots of an
univariate polynomial: we have to minimize the function $x \mapsto
\mu(Q_{\mathsf{E}^*_6} + x R_1)^d / \det(Q_{\mathsf{E}^*_6} + x R_1)$
where we know, because of the approximate solution, that
$\mu(Q_{\mathsf{E}^*_6} + x R_1)$ is a polynomial for all points in a
sufficiently small neighborhood of the exact $x$. This leads to the Ansatz
\[
Q^{c2}_6 = Q_{\mathsf{E}^*_6} + \frac{\sqrt{1057} - 1}{88} R_1.
\]
Now we can use the tools we introduced in Section~\ref{sec:classical}
to prove rigorously that this Ansatz works.

\begin{theorem}
\label{th:secondcovering}
The PQF $Q^{c2}_6$ gives a locally optimal lattice covering. Its
covering density is
\[
\Theta(Q^{c2}_6) = 
\frac{\sqrt{1124895337698\sqrt{1057}-33713139497730}}
{3543122}
\cdot \kappa_6
\approx 
2.466121650.
\]
\end{theorem}

\begin{proof}
Our proof is again computational. We provide the MAGMA program
\texttt{check\_q6c2.m} at the \texttt{arXiv.org} e-print archive.
To access it, download the source files for the paper \texttt{math.MG/0403272}.
Since it uses rational arithmetic only, the proof is rigorous.

Let us describe the steps.  First we compute the Delone subdivision of
$Q^{c2}_6$: We compute the Delone subdivisions of the three PQFs $Q_1
= \frac{9}{10}Q_{\mathsf{E}^*_6} + \frac{3}{10}R_1$, $Q_2 =
\frac{11}{10}Q_{\mathsf{E}^*_6} + \frac{3}{10}R_1$, and $Q_3 =
\frac{11}{10}Q_{\mathsf{E}^*_6} + 
\frac{5}{10}R_1$. 
Then we check that
they coincide and that it turns out to be a triangulation. Then we
show that $Q^{c2}_6 \in \conv\{Q_1,Q_2,Q_3\}$ so that $\Del(Q^{c2}_6)
= \Del(Q_1)$. Now we can compute the circumradii of all the Delone
simplices by the formula~(\ref{eq:circumradius}) in proof of
Proposition~\ref{prop:lmi}. This gives the value~$\Theta(Q^{c2}_6)$.
Finally we prove that $Q^{c2}_6$ gives a locally optimal lattice
covering using the criterion of
Proposition~\ref{prop:localoptcrit}. We compute the gradients $g_L$
for all simplices $L \in \Del(Q^{c2}_6)$ achieving maximum
circumradius. Summing them up yields a multiple of $-(Q^{c2}_6)^{-1}$.
\end{proof}

This computational proof takes about one minute on a $2$GHz Intel
Pentium computer.

What is the general pattern behind the beautification process?  Let
$Q$ be a locally optimal lattice covering with Delone
triangulation~$\MD$.  We use the symmetry of~$\MD$ to find the
subspace in which~$Q$ lies. This reduces the number of involved
variables. The simplices of the Delone triangulation which have
circumradius~$1$ give equality constraints. Then we maximize the
determinant of the quadratic forms lying in the subspace subject to
the equality constraints. For this optimization problem, which
involves only algebraic equations, we can use Gr\"obner basis
techniques.

Unfortunately, we were not able to solve the corresponding algebraic
equations for $Q^{c1}_6$ so we have to be satisfied with an
approximation.

\begin{theorem}
\label{th:firstcovering}
The covering density of the PQF $Q^{c1}_6$ is at most $2.464801$.
\end{theorem}

\begin{proof}
The covering density of the PQF $Q_{\mathsf{E}^*_6} + x R_2 + y R_3 +
z R_4$ with $x = 0.15266328480099$, $y = 0.32884740614948$, $z =
0.13827491241153$ is smaller than $2.464801$. For a computational
proof of this fact we provide the MAGMA program \texttt{check\_q6c1.m},
available from the source files of the paper \texttt{math.MG/0403272}
at the \texttt{arXiv.org} e-print archive.
The steps are similar to the first steps of the proof
of Theorem~\ref{th:secondcovering}.
\end{proof}

On the basis of Theorem~\ref{th:e6d} and our extensive computational
experiments we make the following conjecture.

\begin{conjecture}
The PQF $Q^{c1}_6$ provides the unique least dense lattice covering
in dimension $6$.
\end{conjecture}

Using a similar and more successful beautification process for
$Q^{pc}_6$, we make the Ansatz
\[
Q^{pc}_6 = Q_{\mathsf{E}^*_6} + \frac{\sqrt{798} - 18}{79} R_1
.
\]

\begin{theorem}
\label{th:firstpackingcoveing}
$Q^{pc}_6$ is a locally optimal solution to the lattice
packing-covering problem, lying in the interior of its secondary
cone. Its lattice packing covering constant is
\[
\gamma(Q^{pc}_6) = 
2\sqrt{2\sqrt{798}-56}
\approx
1.411081242
.
\]
\end{theorem}

\begin{proof}
This is similar to the proof of Theorem~\ref{th:secondcovering}. We
provide the MAGMA program \texttt{check\_q6pc.m},
available from the source files of the paper \texttt{math.MG/0403272}
at the \texttt{arXiv.org} e-print archive.
\end{proof}

\subsection{Dimension $7$}
\label{ssec:dim7} 

After analyzing the $6$-dimensional case, we got a feeling of where we
have to search for good $7$-dimensional lattice coverings. We took
$Q_{\mathsf{E}_7^*}$ and a lattice vector which is a longest vector of
the shortest vectors in the cosets $\Z^d/2\Z^d$. We computed the
stabilizer group of this vector and the invariant subspace of this
group. By perturbing $Q_{\mathsf{E}_7^*}$ in this subspace randomly,
we found a PQF whose Delone subdivision is a triangulation. We solved
the determinant maximization problem of Section~\ref{ssec:solvinglcp}
which belongs to this Delone triangulation and found the PQF
\[
Q_7^c =
\begin{pmatrix}
12 &  1 &  1 &  1 &  1 &  1 &  5\\
 1 & 12 &  1 &  1 &  1 &  1 &  5\\
 1 &  1 & 12 &  1 &  1 &  1 &  5\\
 1 &  1 &  1 & 12 &  1 &  1 &  5\\
 1 &  1 &  1 &  1 & 12 &  1 &  5\\
 1 &  1 &  1 &  1 &  1 & 12 & -6\\
 5 &  5 &  5 &  5 &  5 & -6 & 14
\end{pmatrix}.
\]
We are quite surprised that this PQF has rational entries.

\begin{theorem}
$Q^{c}_7$ is a locally optimal solution to the lattice covering
problem, lying in the interior of its secondary cone. Its
inhomogeneous minimum is $\mu = \frac{15}{2}$, its determinant is $\det
Q^{c}_7 = 2 \cdot 11^6$, so that $\Theta(Q_7^c) \approx 2.900024$.
\end{theorem}

\begin{proof}
Again, this is similar to the proof of
Theorem~\ref{th:secondcovering}. We provide the MAGMA program
\texttt{check\_q7c.m}, available from the source files of the paper
\texttt{math.MG/0403272} at the \texttt{arXiv.org} e-print archive.
\end{proof}

\section*{Acknowledgment}
We wish to thank Mathieu Dutour for the simplified formulation of the
optimization problem in Section~\ref{ssec:solvingslpc}, J\"org Rambau
for pointing out the reference to Proposition~\ref{prop:combtria},
Francisco Santos for helpful discussions, and 
Tyrrell B. McAllister for improving the grammar of our text.
In particular we wish to thank one of the
anonymous referees for his detailed report with many helpful comments
and suggestions on a previous version.

\bibliographystyle{amsalpha}
\bibliography{biblio}

\vspace{2ex}

\begin{samepage}
\noindent
\textit{Achill Sch\"urmann, Department of Mathematics, University of
Magdeburg, 39106 Magdeburg, GERMANY, email:
\texttt{achill@math.uni-magdeburg.de}}
\end{samepage}

\vspace{1ex}

\begin{samepage}
\noindent
\textit{Frank Vallentin, Einstein Institute of Mathematics, The Hebrew
University of Jerusalem 91904, ISRAEL, email:
\texttt{frank.vallentin@gmail.com}}
\end{samepage}

\end{document}